\documentclass[11pt]{article}
\usepackage{amssymb,amsfonts,amsmath,amsthm,cite}
\usepackage{graphicx,epsfig,subfigure}
\usepackage{enumerate}
\usepackage{extarrows}
\usepackage{supertabular}
\usepackage{array}
\usepackage{float}

\newtheorem{thm}{Theorem}[section]

\newtheorem{lem}{Lemma}[section]

\makeatletter \@addtoreset{equation}{section}
\makeatletter \@addtoreset{figure}{section}
\makeatletter \@addtoreset{table}{section}

\def\pf{\noindent {\it Proof.\ }}
\def\qed{\hfill \rule{4pt}{7pt}}
\def\ld{{\rm{ldeg}}}
\def\rd{{\rm{rdeg}}}

\begin{document}
\begin{center}
{\Large\bf Zigzag Stacks and $m$-Regular Linear Stacks
} \\[20pt]
{\large William Y.C. Chen$^{a,b,*}$, Qiang-Hui Guo$^a$, Lisa H. Sun$^a$ and Jian Wang$^a$} \\[10pt]
$^{a}$Center for Combinatorics, LPMC-TJKLC \\
Nankai University, Tianjin 300071, P. R. China \\[10pt]
$^{b}$Center for Applied Mathematics\\
Tianjin University, Tianjin 300072, P. R. China\\[10pt]
{Emails: $^*$chen@nankai.edu.cn, guo@nankai.edu.cn\\
sunhui@nankai.edu.cn, wj121313@126.com}
\end{center}

\noindent {\bf Abstract.}  The contact map of a protein fold is a graph that represents the patterns of contacts in the fold.
It is known that the  contact map
 can be decomposed into stacks and queues.
RNA secondary structures are special stacks in which the degree of each vertex is  at most one and each arc has length at least two.
Waterman and Smith derived a formula for the number of  RNA secondary structures of length $n$ with exactly $k$ arcs. H\"{o}ner zu Siederdissen et al. developed a folding algorithm for extended RNA secondary structures in which each vertex has maximum degree two.
An equation for the generating function of extended RNA secondary
structures was obtained by
M\"{u}ller and Nebel by using a  context-free grammar approach, which leads to an asymptotic formula.  In this paper, we consider $m$-regular linear stacks,
where  each arc has length at least $m$ and the degree of
each vertex is bounded by two. Extended RNA secondary structures are exactly $2$-regular linear stacks.  For any $m\geq 2$, we obtain an equation for the generating function
of the $m$-regular linear stacks. For given $m$, we can deduce
a recurrence relation and an asymptotic formula for the number of $m$-regular linear stacks on $n$ vertices.
To establish the equation, we use the reduction
operation of Chen, Deng and Du to transform
 an $m$-regular linear stack to an $m$-reduced  zigzag (or alternating) stack.
Then we find an equation for $m$-reduced zigzag stacks
leading to an equation for $m$-regular linear stacks.

\section{Introduction}
\label{sec-intro}

Proteins are  polymer chains consisting of amino acid residues of twenty  types. The function of a protein is directly dependent on its three dimensional structure. Due to the   complexity of the full-atom protein model, lattice models have been proposed  and
extensively studied. Lattice models often preserve important features of the protein structure, and enable us to focus on dominant aspects of a protein structure. In such a model, protein folds are
 represented by self-avoiding walks on the specific lattice.

\begin{figure}[htp]
  \begin{center}
    \includegraphics[scale=0.8]{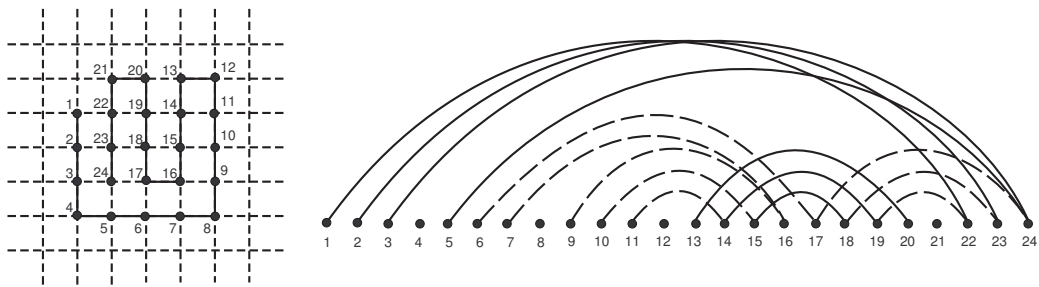}
  \end{center}
  \caption{A protein fold on the 2D square lattice and its contact map.}\label{fig-saw2cm}
\end{figure}

When two amino acids in a protein fold come very close to each other, say, closer than a predetermined threshold,   they presumably form some kind of bond, which is called a contact.
Let $S=s_1s_2\ldots s_n$ represent the amino acid residue sequence of a protein. When we consider the protein fold as a self-avoiding walk on some regular lattice, two residues $s_i$ and $s_j$ are in a contact if they reside on two adjacent points in the lattice,
but not consecutive in the sequence. Let the vertex $i$ stand for the residue $s_i$.
The contact map of a folding of $S$ is a diagram with
vertices $1, 2, \ldots, n$  arranged on a
horizontal line and there is an edge between two vertices
 if they are in contact. See Figure \ref{fig-saw2cm} for an illustration.

 Contacts play a fundamental role in the HP-model for protein folding, see \cite{Dill90,DBYFYTC95}. Contact maps of protein folds have been extensively studied from various perspectives, such as protein folding prediction \cite{Doma00, VKD97}, structure alignment \cite{GIP99,LCW01,AMW07}, protein secondary structure \cite{KL07, Zuk84}, and protein structure data mining \cite{HSSBZ02}.
 Crippen \cite{Crip00} has studied the enumeration of contact maps. Vendruscolo et al. \cite{VKDLS99} investigated  statistical properties of contact maps. Goldman et al. \cite{GIP99} discovered that contact maps for protein folds in two dimension can be decomposed into  ``simpler'' graphs, called stacks and queues.
In combinatorial words, a \emph{stack} is  a noncrossing diagram, and a \emph{queue} is a nonnesting diagram.

\begin{thm}[Goldman et al. \cite{GIP99}]
For any protein sequence $S$, the contact map of any two-dimensional fold of $S$ can be decomposed into (at most) two stacks and one queue.
\end{thm}

For example, the contact map in Figure \ref{fig-saw2cm} can be decomposed into two stacks: $\big\{(6,17)$, $(7,16)$, $(9,16)$, $(10,15)$, $(11,14)$, $(17,24)$, $(18,23)$, $(19,22)\big\}$ and $\big\{(13,20)$, $(14,19)$, $(15,18)\big\}$, and one queue  $\big\{(1,22)$, $(2,23)$, $(3,24)$, $(5,24)\big\}$.

Recently, Agarwal et al. \cite{AMW07} found a  similar decompositions of contact maps of protein folds in the
 three dimensional cubic lattice.

As pointed out by Istrail and Lam,  the enumeration of stacks and of queues is related to an approximation algorithm for computing the partition function of self-avoiding walks in two dimensions. Denote the numbers of stacks on $n$ vertices by $s(n)$. We notice that from the combinatorial intepretation of Schr\"{o}der numbers $a_n$ in terms of  noncrossing graphs, see  \cite[Exercise 6.39(p)]{Stan99},  it is easy to see  that $s(n) = 2^{n-1}a_{n-2}$,
where $a_n$ is the Schr\"{o}der number whose generating function is given by
$$
\sum_{n\geq 0} a_nx^n = \frac{1-x-\sqrt{1-6x+x^2}}{2x}.
$$

When folding a protein on a specific lattice, the lattice model
 leads to  degree and arc length constraints to the corresponding contact map.
For instance, in a contact map on 2D square lattice, each internal vertex has maximum degree $2$ and each arc has length at least $3$. In the case of 2D triangular lattice, each vertex has maximum
 degree $4$ and each arc has length at least $2$. In  the hexagonal lattice, the degree of each vertex is at most 1, but the
  length of each arc is at least $5$. For other lattice models,
  see \cite{PDT08}.  A stack with arc length at least $m$ is called an $m$-regular stack and a stack with degree of each vertex bounded by two is called linear.

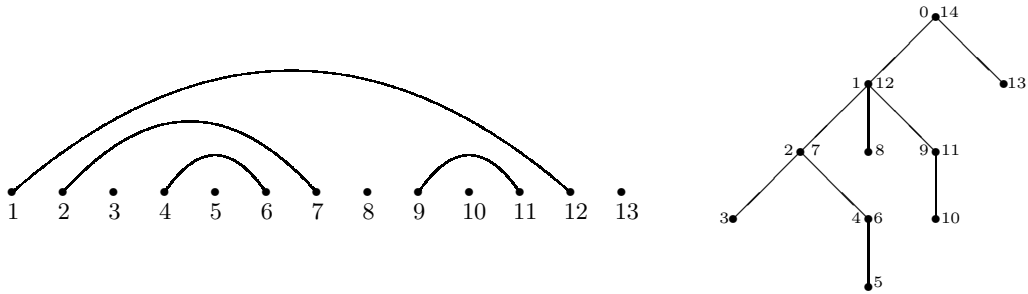
\begin{figure}[ht]
\begin{center}

\setlength{\unitlength}{0.45mm}

\begin{picture}(300,100)

\multiput(2,28)(15,0){13}{\circle*{2}}

\put(0.5,20){\footnotesize $1$}     \put(15.5,20){\footnotesize $2$}

\put(30.5,20){\footnotesize $3$}    \put(45.5,20){\footnotesize $4$}

\put(60.5,20){\footnotesize $5$}    \put(75.5,20){\footnotesize $6$}

\put(90.5,20){\footnotesize $7$}    \put(105.5,20){\footnotesize $8$}

\put(120.5,20){\footnotesize $9$}   \put(135,20){\footnotesize $10$}

\put(150,20){\footnotesize $11$}  \put(165,20){\footnotesize $12$}

\put(180,20){\footnotesize $13$}

\qbezier[1000](2,28)(84.5,100)(167,28) \qbezier[1000](17,28)(54.5,70)(92,28)

\qbezier[1000](47,28)(62,50)(77,28)

\qbezier[1000](122,28)(137,50)(152,28)

\multiput(275,80)(20,0){1}{\circle*{2}}
\multiput(255,60)(40,0){2}{\circle*{2}}
\multiput(235,40)(20,0){3}{\circle*{2}}
\multiput(215,20)(40,0){2}{\circle*{2}}
\multiput(275,20)(20,0){1}{\circle*{2}}
\multiput(255,0)(40,0){1}{\circle*{2}}

\put(275,80){\line(-1,-1){20}}
\put(275,80){\line(1,-1){20}}
\put(255,60){\line(-1,-1){20}}
\put(255,60){\line(1,-1){20}}
\put(255,60){\line(0,-1){20}}
\put(235,40){\line(-1,-1){20}}
\put(235,40){\line(1,-1){20}}
\put(275,40){\line(0,-1){20}}
\put(255,20){\line(0,-1){20}}

\put(270,80){\tiny $0$} \put(276,80){\tiny $14$}
\put(250,59){\tiny $1$} \put(257,59){\tiny $12$}
\put(296,59){\tiny $13$}
\put(230,39){\tiny $2$} \put(238,39){\tiny $7$}
\put(257,39){\tiny $8$}
\put(270,39){\tiny $9$} \put(276.5,39){\tiny $11$}
\put(211,19){\tiny $3$}
\put(250,19){\tiny $4$} \put(256.5,19){\tiny $6$}
\put(276.5,19){\tiny $10$}
\put(256.5,0){\tiny $5$}
\end{picture}
\end{center}
\caption{An RNA secondary structure and its Schmitt-Waterman contact tree.}
\label{secrna}
\end{figure}

   RNA secondary structures can be viewed as $2$-regular stacks with maximum degree $1$. By establishing a bijection between RNA secondary structures and linear trees, Schmitt and Waterman \cite{WS78,SW94} provided an explicit formula for the number of RNA secondary structures on $n$ vertices and $k$ arcs. See  Figure \ref{secrna} for an example of the bijection.

\begin{thm}[Schmitt and Waterman \cite{SW94}]
The number of RNA secondary structures of length $n$ with $k$ arcs is given by
\begin{equation}\label{eqn-RNAnum}
s_2(n,k)=\frac{1}{k}\binom{n-k}{k+1}\binom{n-k-1}{k-1}.
\end{equation}
\end{thm}

In a recent survey   \cite{IL09}, the authors raised the question concerning generalizations of the Schmitt-Waterman counting formulas to stacks  and to queues. In fact, Nebel \cite{Nel02}
derived the generating function  of $m$-regular RNA secondary structures by using the binary trees and the
Horton-Strahler number. Based on the bijection between matchings and walks inside Weyl-chambers given by  \cite{CDDS07}, Jin, Qin and Reidys \cite{JQR08} derived a  formula
 for the number of $k$-noncrossing RNA structures with pseudoknots of length $n$ with $l$ isolated vertices.

\begin{thm} Let
$f_k(n,l)$ be the number
of $k$-noncrossing digraphs over $n$ vertices with exactly $l$ isolated vertices. Then the number of RNA
structures of pseudoknot type $k-2$ with $l$ isolated vertices is given by
\begin{equation}\label{RNA-1}
s_k(n,l)=\sum_{b=0}^{(n-l)/2}(-1)^b{n-b \choose b}f_k(n-2b,l),
\end{equation}
where $n\geq 1$, $1\leq l \leq n$ and $k\geq 2$.
\end{thm}

In particular, when $k=2$, formula \eqref{RNA-1} reduces to  Schmitt and Waterman's formula \eqref{eqn-RNAnum}. When $k=3$, Jin and Reidys \cite{JR07} obtained the following asymptotic formula for the $3$-noncrossing RNA structures
\[
S_3(n)\thicksim \frac{10.4724\cdot 4!}{n(n-1)\cdots (n-4)} \bigg(\frac{5+\sqrt{21}}{2}\bigg)^n.
\]

In 2011, H\"{o}ner zu Siederdissen et al. \cite{HBSH11} presented a model of extended RNA secondary structures in which  the degree of each vertex is bounded by two. They provided a folding algorithm which is known to be the first thermodynamics-based algorithm that allows the degrees of vertices to be two.
M\"{u}ller and Nebel \cite{MN13} studied the enumeration of extended RNA secondary structures by using a context-free grammar approach. They obtained an equation satisfied by the ordinary generating functions of the number of extended RNA secondary structures.
\begin{thm} Let
\[
S(z)=\sum_{n\ge 0} r_2(n)z^n,
 \]
 where $r_2(n)$ is the number of extended RNA secondary structures of length $n$ and $r_2(0)=0$. Then we have
\begin{align}\label{gf-exrna}
S(z) = &\ 4z^5S^5(z)+ (4z^3-7z^4+9z^5)S^4(z)+(-8z^2+11z^3-14z^4+7z^5)S^3(z)\nonumber\\[5pt]
&\quad  +(5z-10z^2+14z^3-9z^4+2z^5)S^2(z)+(3z-7z^2+7z^3-2z^4)S(z)\nonumber\\[5pt]
&\quad +z-2z^2+z^3
\end{align}
and
 \begin{equation}\label{asyr2}
 r_2(n)\thicksim 0.250536155\times 4.1012475^n \cdot n^{-\frac{3}{2}}.
 \end{equation}
 \end{thm}

 In this paper, we are mainly concerned with the $m$-regular linear stacks, in which the arc length is at least $m$ with $m\ge 2$ and the degree of each vertex is bounded by $2$. For example, the two stacks in Figure \ref{fig-saw2cm} are $3$-regular linear stacks.
 Extended RNA secondary structures are exactly $2$-regular linear stacks.

 We obtain an equation for the generating function of the $m$-regular linear stacks. For given $m$, we can derive an explicit recurrence relation and an asymptotic formula  for the number of $m$-regular linear stacks of length $n$.

\begin{figure}[h]
\begin{center}
\setlength{\unitlength}{0.4mm}
\begin{picture}(120,60)
\multiput(0,10)(20,0){7}{\circle*{2}}
\qbezier[1000](0,10)(60,80)(120,10)
\qbezier[1000](20,10)(70,65)(120,10)
\qbezier[1000](20,10)(50,45)(80,10)
\qbezier[1000](40,10)(60,35)(80,10)
\put(-2,0){\small $1$}\put(18,0){\small $2$}\put(38,0){\small $3$}
\put(58,0){\small $4$}\put(78,0){\small $5$}\put(98,0){\small $6$}
\put(118,0){\small $7$}
\end{picture}
\end{center}
\caption{A zigzag stack.}
\label{fig-alterStack}
\end{figure}
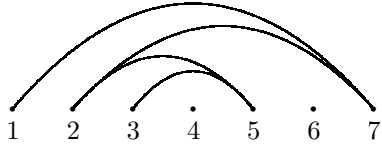

Using the reduction operation on the standard representation of
an $m$-regular partition due to Chen et al. \cite{CDD05}, we find a  class of zigzag stacks  which are  in one-to-one correspondence with  $m$-regular linear stacks.
More precisely, a zigzag stack, or an alternating stack, is  a noncrossing diagram
without multiple edges or loops subject to the following conditions:
\begin{enumerate}
\item[(1)] The degree of each vertex is bounded by 2;

\item[(2)] For each vertex $v$ of degree 2, the two arcs are on the same side with respect to the position of $v$.
\end{enumerate}
 Notice that isolated vertices are allowed in a
 zigzag stack.
 Figure \ref{fig-alterStack} gives a zigzag stack.
It is easy to see that each connected component of a zigzag stack forms a zigzag path.

\begin{figure}[ht]
\begin{center}
\setlength{\unitlength}{0.5mm}
\begin{picture}(240,40)
\multiput(0,10)(12,0){8}{\circle*{2}}
\multiput(160,10)(16,0){6}{\circle*{2}}
\qbezier[1000](0,10)(24,48)(48,10)
\qbezier[1000](12,10)(30,36)(48,10)
\qbezier[1000](48,10)(66,36)(84,10)
\put(-2,0){$1$}\put(10,0){$2$}\put(22,0){$3$}
\put(34,0){$4$}\put(46,0){$5$}\put(58,0){$6$}
\put(70,0){$7$}\put(82,0){$8$}
\put(102,10){$\xlongrightarrow{\ \rm reduction\  }$}
\qbezier[1000](160,10)(176,36)(192,10)
\qbezier[1000](176,10)(184,24)(192,10)
\qbezier[1000](224,10)(232,24)(240,10)
\put(158,0){$1$}\put(174,0){$2$}\put(190,0){$3$}
\put(206,0){$4$}\put(222,0){$5$}\put(238,0){$6$}
\end{picture}
\end{center}
\caption{The reduction from a non-linear stack to a zigzag stack.}
\label{fig-reduce}
\end{figure}
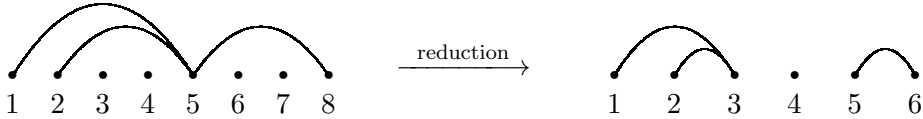

Given a vertex $v$ in a stack, we denote the degree of $v$  by $\deg(v)$, and denote the left-degree and right-degree of $v$ by $\ld(v)$ and $\rd(v)$, respectively.
We find  that the  reduction operation transforms an $m$-regular linear stack of length $n+m-1$ to a zigzag stack of length $n$ with the following two constraints:
\begin{enumerate}
  \item[(1)] For  $1\le i\le n-m+1$, $\ld(i)+\rd(i+m-1)\leq 2$;
  \item[(2)] For  $1\le i<j\le n$, if $\ld(i)>0$ and  $\rd(j)>0$, then $j-i\geq m-1$.
\end{enumerate}
A zigzag stack satisfying the above two constraints is
 called an $m$-reduced zigzag stack. The above conditions
can be used to characterize the substructures of $m$-reduced zigzag stacks. For example, the zigzag stack in Figure \ref{fig-reduce} is obtained from a non-linear $3$-regular stack  by applying the reduction
operation twice.
Thus any $3$-reduced zigzag stack can not contain a
substructure like the reduced zigzag stack in Figure \ref{fig-reduce}. We shall use
Conditions (1) and (2) to describe the substructures in the
decomposition of $m$-reduced zigzag stacks. As will be seen,
the substructure in Figure \ref{fig-reduce} is not a valid substructure according to the characterizations as given in Theorem \ref{thm-subs}.

Denote the set of $m$-regular linear stacks and $m$-reduced zigzag stacks of length $n$ by $\mathcal{R}_m(n)$ and $\mathcal{Z}_m(n)$, respectively. Clearly, when $m=2$, $\mathcal{R}_2(n)$ induces to the extended RNA secondary structures.
Denote $r_m(n)=|\mathcal{R}_m(n)|$ and $z_m(n)=|\mathcal{Z}_m(n)|$, and let
 \[
R_m(x)=\sum_{n=0}^\infty r_m(n)x^n,
\]
\[
Z_m(x)=\sum_{n=0}^\infty z_m(n)x^n,
\]
where $r_m(0)=z_m(0)=1$.
Then the reduction algorithm implies that $z_m(n)=r_m(n+m-1)$ and
\begin{equation}\label{re-RmZm}
R_m(x)=1+x+x^2+\cdots+x^{m-2}+x^{m-1}Z_m(x).
\end{equation}

By decomposing an $m$-reduced zigzag stack into a connected component and a list of substructures,  we derive an equation satisfied by $Z_m(x)$. Furthermore, by applying relation \eqref{re-RmZm}, we obtain that the following relation for $R_m(x)$ with polynomial coefficients.

\begin{thm}\label{thm22} We have
\begin{equation}\label{gf-r-l}
c_5(x)R_m^5+c_4(x)R_m^4+c_3(x)R_m^3+c_2(x)R_m^2+c_1(x)R_m+c_0(x)=0,
\end{equation}
where
\begin{align*}
c_0(x)&=(x-1) (x^m-1)^3,\\[5pt]
c_1(x)&= (x^m-1)^2(x^{2m+1}-2x^{m+2}-x^{m+1}+ x^m - 3 x^3 + 8 x^2 - 3 x -1),\\[5pt]
c_2(x)&= -x(x-1) (x^m-1)(5x^{2m+1}-6x^{m+2}
-9x^{m+1}+5x^m-3x^3+12x^2+x-5),\\[5pt]
c_3(x)&= x^2(x-1)^2(x^m-1)(11x^{m+1}-8x^2-11x+8),\\[5pt]
c_4(x)&= x^3(x-1)^3(-11x^{m+1}+4x^2+11 x-4),\\[5pt]
c_5(x)&=4 x^5(x-1)^4.
\end{align*}
\end{thm}

From the above equation \eqref{gf-r-l}, we can deduce the recurrence relation for $r_m(n)$ for given $m$. Applying Newton-Puiseux Expansion Theorem \cite{FS09}, we are led to the
following asymptotic formula
\[
r_m(n)\thicksim \gamma\cdot \omega^n \cdot n^{-\frac{3}{2}},
\]
where $\gamma$ and $\omega$ are constants.
When $m=2$, it coincides with M\"{u}ller and Nebel's formula \eqref{asyr2}  for  extended RNA secondary structures, see  \cite{MN13}. For $m=3,4,5,6$, we have
\begin{align*}
r_3(n) &\thicksim 0.19005341 \times 3.5271506^n \cdot n^{-\frac{3}{2}},\\[5pt]
r_4(n)& \thicksim  0.145636571 \times 3.2431591^n \cdot n^{-\frac{3}{2}},\\[5pt]
r_5(n)&\thicksim 0.112004701\times 3.0833083^n \cdot n^{-\frac{3}{2}},\\[5pt]
r_6(n)&\thicksim 0.086237333\times 2.9880679^n\cdot  n^{-\frac{3}{2}}.
\end{align*}

\section{Zigzag stacks}\label{sec-zigzag}

In this section, we derive an equation satisfied by
the generating function of the number of zigzag stacks
on $n$ vertices.  To enumerate zigzag stacks, we introduce
the primary component decomposition.

\begin{figure}[ht]
\begin{center}
\setlength{\unitlength}{0.4mm}
\begin{picture}(160,70)
\multiput(0,10)(25,0){7}{\circle*{2}}
\multiput(6,8)(25,0){7}{\framebox(13,4)}
\qbezier[1000](0,10)(75,110)(150,10)
\qbezier[1000](0,10)(37.5,60)(75,10)
\qbezier[1000](25,10)(50,44)(75,10)
\qbezier[1000](25,10)(37.5,27)(50,10)
\qbezier[1000](100,10)(125,44)(150,10)
\qbezier[1000](100,10)(112.5,27)(125,10)
\put(-3,0){\small $v_1$}\put(22,0){\small $v_2$}\put(47,0){\small $v_3$}
\put(72,0){\small $v_4$}\put(97,0){\small $v_5$}\put(122,0){\small $v_6$}
\put(147,0){\small $v_7$}
\end{picture}
\end{center}
\caption{Decomposition of a zigzag stack.}
\label{fig-subinterval}
\end{figure}
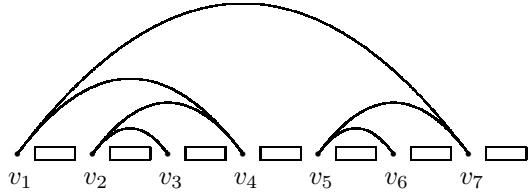

Denote $[n]=\{1, 2, \ldots, n\}$.
Let $S$ be a zigzag stack on $[n]$,
and define the primary component of $S$ to be the
connected component
containing the vertex $1$. The primary component will split $[n]$ into disjoint intervals, on which smaller zigzag stacks can be constructed. This enables us to establish a recursive procedure to enumerate zigzag stacks, so that we can derive an equation on the generating function for the number of zigzag stacks.
Figure \ref{fig-subinterval}  illustrates a primary component decomposition of a zigzag stack.

The following lemma shows that the primary component decomposition
leads to a primary component along with
zigzag stacks on the intervals.

\begin{lem}\label{lem-pc}
Let $S$ be a zigzag stack, and $C$ be the primary component of $S$.
 Then there is no
arc of $S$ that connects two vertices in different intervals.
In other words, a zigzag stack can be decomposed into a primary component along with a list of zigzag stacks on the intervals.
\end{lem}

\pf Clearly, any vertex in an interval is not connected to any
vertex of $C$. Assume to the contrary that the lemma does not
hold. Then there exists an arc $(i,j)$ with $i<j$  such that $i\in I$ and $j\in J$, where $I$ and $J$ are two different intervals of $S$.
  Let $k$ be the vertex of $C$ such that $k$ is next to the
  last vertex of $I$. We now have $1<i<k<j$.  Since $1$ and $k$ are  in the connected component $C$, the arc $(i,j)$ must intersect with some arc of $C$, contradicting the assumption that $S$ is a
  stack. This completes the proof.
\qed

Let $c_n$ denote the number of connected zigzag stacks on  $[n]$.
We have the following formula.

\begin{lem}\label{lem-cas2n} If $n\ge 2$, we have
\begin{equation}\label{eq-cas2n}
c_n=n-1.
\end{equation}
\end{lem}

\begin{figure}[ht]
\begin{center}
\setlength{\unitlength}{0.4mm}
\begin{picture}(230,50)
\put(2,10){\circle*{2}} \put(17,10){\circle*{2}}
\put(32,10){\circle*{2}} \put(47,10){\circle*{2}}
\put(77,10){\circle*{2}} \put(92,10){\circle*{2}}
\put(107,10){\circle*{2}} \put(122,10){\circle*{2}}
\put(162,10){\circle*{2}} \put(177,10){\circle*{2}}
\put(192,10){\circle*{2}} \put(207,10){\circle*{2}}
\put(222,10){\circle*{2}} \put(237,10){\circle*{2}}
\qbezier[1000](2,10)(24,40)(47,10)
\qbezier[1000](17,10)(32,30)(47,10)
\qbezier[1000](17,10)(24,20)(32,10)
\qbezier[1000](77,10)(99,40)(122,10)
\qbezier[1000](77,10)(92,30)(107,10)
\qbezier[1000](92,10)(99,20)(107,10)
\qbezier[1000](162,10)(199,60)(237,10)
\qbezier[1000](162,10)(177,30)(192,10)
\qbezier[1000](177,10)(184,20)(192,10)
\qbezier[1000](207,10)(222,30)(237,10)
\qbezier[1000](207,10)(214,20)(222,10)
\put(20,-4){\small (a)}
\put(95,-4){\small (b)}
\put(195,-4){\small (c)}
\end{picture}
\end{center}
\caption{Three cases of connected zigzag stacks}
\label{threecases-zigzag}
\end{figure}
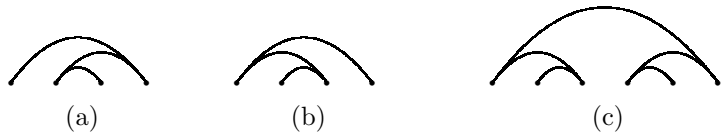

\pf The lemma is obvious true for $n= 2, 3$.
For $n\geq 4$, let $C$ be a connected zigzag stack
on $[n]$. Clearly, $(1,n)$ is an arc in $C$, since
$C$ is connected and is zigzag. There are three cases with respect to degrees of $1$ and $n$, see Figure \ref{threecases-zigzag}.

 If $\deg(1)=1$, as shown in Figure \ref{threecases-zigzag} (a), $C$ is the unique diagram  consisting of the arcs $(1, n)$, $(2,n)$, $(2, n-1)$, and so on. Similarly, if $\deg(n)=1$, $C$ is uniquely determined.

If $\deg(1)=\deg(n)=2$,
assume that $(1, i)$ is an arc of $C$, where $i<n$.
Then clearly, $(i+1, n)$ is an arc of $C$.
Moreover, $C$ is   determined once the arcs $(1,i)$ and $(i+1,n)$\ are given. Since $i$ can be any vertex in $\{2, 3, \ldots, n-2\}$, there are $n-3$ choices of $C$ in this case.
In summary, there are a total number of $n-1$ connected zigzag stacks on $[n]$.   \qed

We set $c_0=c_1=1$.
For $n\geq 1$, let $\mathcal{Z}(n)$ denote
the set of zigzag stacks on $[n]$, and let $z(n)=|\mathcal{Z}(n)|$.
We set
$z(0)=1$. Let
 \[
 Z(x)=\sum_{n\geq 0}z(n)x^n
 \]  denote the ordinary generating function of $z(n)$. The following theorem gives an equation satisfied by $Z(x)$.

\begin{thm}\label{thm-2-1}
We have
\begin{align}
& x^2(x-1)Z^3(x)+2xZ^2(x)-(x+1)Z(x)+1=0. \label{eq-A2alg}
\end{align}
\end{thm}

\pf  Let $S$ be a zigzag stack, and let $C$ be the
primary component of $S$.
Assume that $C$ has $k$ vertices,  then
$S$ can be decomposed into a primary component  $C$ and
$k$ zigzag stacks $S_1, S_2, \ldots, S_k$, where
$S_i$ are allowed to be empty.

Let $d_i$ be the number of vertices in $S_i$. Then we have
 \begin{align}
 z(n)&= \sum_{k=1}^{n}c_k \sum_{d_1+d_2+\cdots+d_k=n-k}
 z(d_1)z(d_2)\cdots z(d_k). \label{eq-s2n}
\end{align}
Multiplying both sides of \eqref{eq-s2n} by $x^n$ and
summing over $n$, we obtain
\allowdisplaybreaks
\begin{align*}
Z(x) &= \sum_{n=0}^\infty \sum_{k=1}^{\infty}c_k \sum_{d_1+d_2+\cdots+d_k=n-k}
 z(d_1)z(d_2)\cdots z(d_k) x^n\\[5pt]
  &=\sum_{k=1}^\infty c_k x^k \sum_{c_1=0}^\infty z(d_1) x^{d_1} \sum_{d_2=0}^\infty z(d_2) x^{d_2}\cdots \sum_{d_k=0}^\infty z(d_k) x^{d_k}\\[5pt]
  &=\sum_{k=1}^\infty c_k \big(xZ(x)\big)^k.
\end{align*}
Since $c_0=c_1=1$ and $c_k=k-1$ for $k\geq 2$, we deduce that
\allowdisplaybreaks
\begin{align*}
Z(x) &= c_0+c_1 xZ+\sum_{k=2}^\infty c_k (xZ(x))^k \\[5pt]
&= 1+xZ(x)+\sum_{k=2}^\infty (k-1) (xZ(x))^k\\[5pt]
&=1+xZ(x)+(xZ(x))^2\sum_{k=0}^\infty (k+1)(xZ(x))^k\\[5pt]
&=1+xZ(x)+\frac{(xZ(x))^2}{(1-xZ(x))^2},
\end{align*}
which yields (\ref{eq-A2alg}).  \qed

Given the   equation \eqref{eq-A2alg}, it is straightforward to derive a second order differential equation for $Z(x)$ and a recurrence relation of $z(n)$, see, for example \cite[Chapter 6]{Stan99}. The computation can also be carried out
by using the Maple package {\tt gfun}, see \cite{SZ94}. Moreover, the asymptotic formula for $z(n)$ can be derived by applying Newton-Puiseux Expansion Theorem \cite{FS09}.

\begin{thm} The differential equation satisfied by $Z(x)$ is as follows
\begin{align*}
{x}^{2} (23&{x}^{3} -26{x}^{2}+ 23
x-4 )( 4{x}^{2}+x-1 )( x-1 )Z''(x)\\[5pt]
&+x (368{x}^{6}-433{x}^{5}+108{x}^{4}+260{x}^{3}-258{x}^
{2}+93x-10 ) Z'(x)\\[5pt]
&\qquad\quad+(184{x}^{6}-87{x}^{5}-117{x}^{4}+217{x}^{3}-129{x}^{2}+30
x-2) Z(x)\\[5pt]
&\qquad\quad\qquad\quad-2  (25{x}^{2}-8x+1) (x-1) =0.
\end{align*}
The number of zigzag stacks on $[n]$  satisfies the following recurrence relation
\begin{equation}\label{eq-AS22rec}
\sum_{i=0}^6 p_i(n)z(n+i)=0,
\end{equation}
where
$$
\begin{array}{ll}
p_0(n)=184+276n+92n^2, & p_1(n)=-520-606n-173n^2, \\[5pt]
p_2(n)=347+480n+124n^2, & p_3(n)=937+210n-10n^2, \\[5pt]
p_4(n)=-1881-678n-60n^2, & p_5(n)=1115+372n+31n^2, \\[5pt]
p_6(n)=-182-54n-4n^2. &
\end{array}
$$
We also have
$$
z(n)\thicksim 0.4781905\times 4.6107186^n\cdot n^{-\frac{3}{2}}.
$$
\end{thm}

The first few values of $z(n)$ are given below
\noindent\begin{center}
\begin{tabular}{|c|c|c|c|c|c|c|c|c|c|c|c|c|}
  \hline
  $n$            & 1 & 2 & 3 & 4 & 5  & 6  & 7   & 8   & 9&10 &11 & 12     \\ \hline
  $z(n)$ & 1 & 2 & 6 & 20 & 70 & 255 & 959 & 3696 & 14520 & 57930
& 234080&  955999\\
  \hline
\end{tabular}
\end{center}

\section{Reduction of $m$-Regular linear stacks}\label{sec-enumS2}

In this section, we use the reduction operation in \cite{CDD05} to transform a $m$-regular linear stack  to a
zigzag stack. Then we give a characterization of
zigzag stacks that are in one-to-one correspondence with
$m$-regular linear stacks.  In the next three sections, we obtain
an equation for the generating function of the number of
$m$-reduce zigzag stacks on $[n]$, which leads to an
equation for the generating function of the number
of $m$-regular linear stacks on $[n]$.

 An $m$-reduced zigzag stack in $\mathcal{Z}_m(n)$ is a zigzag stack
  satisfying the following two conditions:
\begin{enumerate}
  \item[(1)] For $1\le i\le n-m+1$, $\ld(i)+\rd(i+m-1)\leq 2$;
  \item[(2)] For  $1\le i<j\le n$, if $\ld(i)>0$ and  $\rd(j)>0$, then $j-i\geq m-1$.
\end{enumerate}

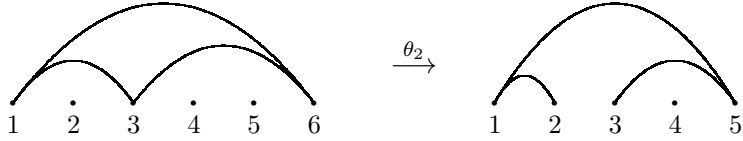
\begin{figure}[ht]\begin{center}\setlength{\unitlength}{0.4mm}
\begin{picture}(250,50)
\put(15,10){\circle*{2}}
\put(35,10){\circle*{2}}
\put(55,10){\circle*{2}}
\put(75,10){\circle*{2}}
\put(95,10){\circle*{2}}
\put(115,10){\circle*{2}}
\put(175,10){\circle*{2}}
\put(195,10){\circle*{2}}
\put(215,10){\circle*{2}}
\put(235,10){\circle*{2}}
\put(255,10){\circle*{2}}
\qbezier[1000](15,10)(35,38)(55,10)
\qbezier[1000](55,10)(85,48)(115,10)
\qbezier[1000](15,10)(65,76)(115,10)
\qbezier[1000](175,10)(185,28)(195,10)
\qbezier[1000](215,10)(235,38)(255,10)
\qbezier[1000](175,10)(215,76)(255,10)
\put(135,20){$\stackrel{\quad \theta_2\quad }{\longrightarrow}$}
\put(13,0){\small $1$}
\put(33,0){\small $2$}
\put(53,0){\small $3$}
\put(73,0){\small $4$}
\put(93,0){\small $5$}
\put(113,0){\small $6$}
\put(173,0){\small $1$}
\put(193,0){\small $2$}
\put(213,0){\small $3$}
\put(233,0){\small $4$}
\put(253,0){\small $5$}
\end{picture}
\end{center}
\caption{An example of the reduction for $m=2$.}\label{fig-reduction}
\end{figure}

The reduction operation $\theta_m$ is defined by
\begin{align} \label{remap}
\theta_m\colon  \mathcal{R}_m(n+m & -1)\longrightarrow \mathcal{Z}_m(n),
\end{align}
in which for  $S\in \mathcal{R}_m(n+m-1)$,  $\theta_m(S)$ is obtained from $S$ by replacing each arc $(i,j)$  by an arc  $(i,j-m+1)$ and   deleting the vertices $n+1,n+2,\ldots,n+m-1$ afterwards.
 Figure \ref{fig-reduction} is  an example for $m=2$.

The following theorem shows that $\theta_m$ is a bijection.

\begin{thm}\label{lem-s2n-as2n}
The map $\theta_m$ is a bijection
between $m$-regular linear stacks on $[n+m-1]$ and $m$-reduced zigzag stacks on $[n]$.
\end{thm}

\pf Let $S$ be any $m$-regular linear stack in $\mathcal{R}_m(n+m-1)$ and $T=\theta_m(S)$. Let $(i,j)$ be an arc of $S$.
Since $S$ is $m$-regular, that is
$j-i\ge m$, we see that the reduced pair $(i, j-m+1)$
is an arc.
  We claim that $T$ is a stack.
Assume to the contrary that $T$ contains two crossing arcs: $(i_1,j_1)$ and $(i_2,j_2)$ with $i_1<i_2<j_1<j_2$.
Thus $(i_1,j_1+m-1)$ and $(i_2,j_2+m-1)$ are arcs of $S$.
Moreover, these two arcs form a crossing since $i_1<i_2<j_1+m-1<j_2+m-1$, a contradiction.
This proves that $T$ is a stack.

Next we show that $T$  is a zigzag stack. We claim that for each vertex $j$ of $T$, either $\ld_T(j)=0$ or $\rd_T(j)=0$. Otherwise, $T$  contains two arcs: $(i,j)$ and $(j,k)$ with $i<j<k$. These
two arcs in $T$ correspond to arcs  $(i,j+m-1)$ and $(j,k+m-1)$,
which form a crossing, since $i<j<j+m-1<k+m-1$. Thus the claim
is proved.

To prove that the degree of each vertex in $T$ is bounded by 2,
we notice that for any vertex $i$ in $T$,
\begin{equation}\label{ldeg}
\ld_T(i)=\ld_S(i+m-1)
\end{equation} and
\begin{equation}\label{rdeg}
\rd_T(i)=\rd_S(i).
\end{equation}
Moreover, if $\ld_T(i)=0$, then $\deg_T(i)=\rd_S (i)$, and if $\rd_T(i)=0$, then $\deg_T(i)=\ld_S (i+m-1)$.
In view of the above claim that for any vertex $i$ in $T$,   $\rd_T(i)=0$ or $\rd_T(i)=0$, we deduce that either $\deg_T(i)=\ld_S (i+m-1)$ or $\deg_T(i)=\rd_S(i)$.
 Given that $S$ is linear, both  $\ld_S (i+m-1)$ and $\rd_S(i)$ are bounded by $2$.
 It follows that $\deg_T(i)\le 2$. Hence $T$ is a zigzag stack.

Employing relations \eqref{ldeg} and \eqref{rdeg}, we deduce that
\begin{align*}
\ld_T(i)+\rd_T(i+m-1)&=\ld_S(i+m-1)+\rd_S(i+m-1)\\[5pt]
&=\deg_S(i+m-1)\leq 2,
\end{align*}
which is Condition (1)  in the definition of an $m$-reduced zigzag stack.

To verify Condition (2), we assume to the contrary that there exist two vertices $j,k$ in $T$ with $j<k$  such that
$\ld_T(j)>0$ and $\rd_T(k)>0$, but $k-j< m-1$. Since $\ld_T(j)>0$ and $\rd_T(k)>0$, there   exist two arcs $(i,j)$, $(k,l)$ in $T$  which are obtained from arcs $(i,j+m-1)$ and $(k,l+m-1)$ in $S$.   Since $k-j<m-1$, we have $i<k<j+m-1<l+m-1$, so that we get two
crossing arcs in $S$, a contradiction. This implies that $T$ satisfies  Condition (2) in the definition of an $m$-reduced zigzag stack. In summary, we have shown that $\theta_m$ is well-defined.

Let $T\in\mathcal{Z}_m(n)$.
To prove that $\theta_m$ is a bijection, we define the inverse map
\[
\varphi_m\colon \ \mathcal{Z}_m(n)\longrightarrow\mathcal{R}_m(n+m -1)
\]
by $\varphi_m(T)=S$, where $S$ is a diagram on $[n+m-1]$ whose arcs are obtained by expanding each arc $(i,j)$ of $T$ to an arc $(i,j+m-1)$. Since $j-i\geq 1$ for each arc $(i,j)$ in $T$, we have $j+m-1-i\geq m$, which says that $S$ is $m$-regular.

Next, to prove $S$ is linear, we notice that for any vertex $i$ in $S$,
  \begin{align}
&\ld_S(i)=\ld_T(i-m+1),\label{ldegs}\\[5pt]
&\rd_S(i)=\rd_T(i).\label{rdegs}
\end{align}
When $1\leq i \leq  m-1$, it is easy to see that $\ld_S(i)=0$ since $S$ is $m$-regular. Then $\deg_S(i)=\rd_T(i)\leq 2$. When $m\le i\le n$, by using Condition (1), we have that
\[
\deg_S(i)=\ld_S(i)+\rd_S(i)=\ld_T(i-m+1)+\rd_T(i)\leq 2.
\]
When $n+1\leq i\leq n+m-1$, it is clear that $\rd_S(i)=0$, which implies that $\deg_S(i)=\ld_T(i-m+1)\leq 2$.  Thus the degree of each vertex in $S$ is bounded by $2$, namely, $S$ is linear.

To prove $S$ is an $m$-regular linear stack, we still need
to show that $S$ is a stack. Otherwise, suppose that $S$ contains two crossing arcs $(i_1,j_1)$ and $(i_2,j_2)$ with $i_1<i_2<j_1<j_2$.
Thus we get two arcs $(i_1,j_1-m+1)$ and $(i_2,j_2-m+1)$ in $T$. If $j_1-m+1>i_2$, then $(i_1,j_1-m+1)$ and $(i_2,j_2-m+1)$  form a crossing in $T$, contradicting
the assumption that $T$ is a stack. If $j_1-m+1=i_2$, then for  vertex $i_2$, we have $\ld_T(i_2)>0$ and $\rd_T(i_2)>0$, contradicting the assumption that  $T$ is a zigzag stack.  If $j_1-m+1<i_2$, then we have $\ld_T(j_1-m+1)>0, \rd_T(i_2)>0$ and $i_2-(j_1-m+1)<j_1-(j_1-m+1)=m-1$.  But this violates  Condition (2) for $T$. So we conclude that by no means does $S$ contain a crossing, namely, $S$ is a stack.

We now have shown that $S$ is an $m$-regular linear stack, that is, $\varphi_m$ is well-defined. It is easily
 seen that $\varphi_m(\theta_m(S))=S$.
 Hence $\theta_m$ is a bijection and the proof is complete.
 \qed

\section{Decomposition of an $m$-reduced zigzag stack}\label{sec-mred}

Based on the bijection given in Theorem \ref{lem-s2n-as2n}, we can transform an $m$-regular linear stack to an $m$-reduced zigzag stack. The advantage of enumerating $m$-reduced
zigzag stack lies in that there is no restriction on the arc lengths. We apply
the primary component decomposition to enumerate $m$-reduced zigzag stacks.

Let
$u$ and $v$ $(u<v)$ be two adjacent vertices in the
primary connected component of $S$. We use
 $\langle u,v\rangle$ to denote the interval of integers
 between $u$ and $v$, that is, $\{u+1, u+2, \ldots, v-1\}$. Note that $\langle u,v\rangle$ is allowed to be empty.
 We use this notation $\langle u,v\rangle$ to distinguish with
  the notation $(u,v)$ for an arc.
If the primary component $C$  does not
 include the last vertex $n$ of $S$, then the
  vertices after the last vertex of $C$
  also form an interval. Since the intervals we are concerned with
   are determined by the primary component $C$, we call such intervals $C$-intervals. Note that the $C$-intervals
   are allowed to be empty. However, if $C$ contains the last vertex $n$,
   we do not consider the empty set after the vertex $n$ as an interval. In this case, if $C$ contains $k$ vertices, then
   there are $k-1$ $C$-intervals.

\begin{figure}[ht]
\begin{center}
\setlength{\unitlength}{0.4mm}
\begin{picture}(250,80)

\multiput(0,10)(16,0){8}{\circle*{2}}
\qbezier[1000](16,10)(24,22)(32,10)
\qbezier[1000](80,10)(72,22)(64,10)
\qbezier[1000](126,10)(158,45)(190,10)
\qbezier[1000](142,10)(150,22)(158,10)
\qbezier[1000](0,10)(95,105)(190,10)
\qbezier[1000](0,10)(48,60)(96,10)
\qbezier[1000](238,10)(246,22)(254,10)
\qbezier[1000](142,10)(158,32)(174,10)
\put(-2,0){\small $1$}\put(14,0){\small $2$}\put(30,0){\small $3$}
\put(46,0){\small $4$}\put(62,0){\small $5$}\put(78,0){\small $6$}
\put(94,0){\small $7$}\put(110,0){\small $8$}

\multiput(126,10)(16,0){9}{\circle*{2}}
\put(124,0){\small $9$}\put(139,0){\small $10$}\put(155,0){\small $11$}
\put(171,0){\small $12$}\put(187,0){\small $13$}\put(203,0){\small $14$}\put(219,0){\small $15$}\put(235,0){\small $16$}\put(251,0){\small $17$}

\end{picture}
\end{center}
\caption{The primary component decomposition of a $3$-reduced zigzag stack.}
\label{fig-primarydecom}
\end{figure}
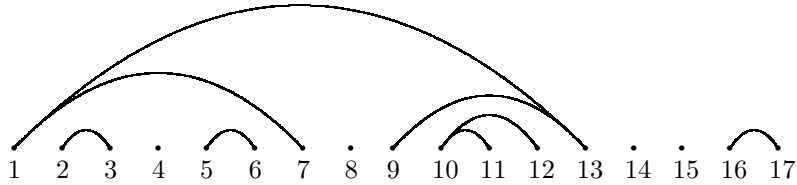

   Figure \ref{fig-primarydecom} illustrates a primary component decomposition of a $3$-reduced zigzag stack $S$. The primary component $C$ is a connected zigzag stack on $\{1,7,9,13\}$, it decomposes  $S$ into four intervals $I_1=\{2,3,4,5,6\}$, $I_2=\{8\}$, $I_3=\{10,11,12\}$ and the last interval is $I_4=\{14,15,16,17\}$.

Let $S$ be a zigzag stack on $[n]$, and let $C$ be the primary component of $S$. Let $I_1, I_2, \ldots, I_k$ be
the $C$-intervals of $S$. Define the sets $J_1, J_2,
\ldots, J_k$ as follows. There are two cases. Case 1: $C$
contains the vertex $n$. Assume that for $1\leq t\leq k$, $I_t$ is the
 interval between the vertices $u_t$ and $v_t$ in $C$.
 Then let $J_t=\{u_t, u_t+1, \ldots, v_t\}$ for each $t$.
 Case 2: $C$ does not contain the vertex $n$.
 Assume that for $1\leq t\leq k-1$,
  $I_t$ is the interval between the vertices $u_t$ and $v_t$ in $C$ and that $I_k$ is the interval after the vertex $u_k$.
   Then let $J_t=\{u_t,u_t+1, \ldots, v_t\}$ for $1\leq t\leq k-1$ and let $J_k=\{u_k,u_k+1, \ldots, n\}$. For example, in Figure \ref{fig-primarydecom}, we see that $J_1=\{1,2,3,4,5,6,7\}$, $J_2=\{7,8,9\}$, $J_3=\{9,10,11,12,13\}$ and $J_4=\{13,14,15,16,17\}$.

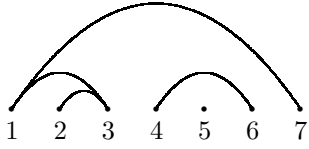
\begin{figure}[ht]
\begin{center}
\setlength{\unitlength}{0.4mm}
\begin{picture}(100,60)

\multiput(0,10)(16,0){7}{\circle*{2}}

\qbezier[1000](48,10)(64,34)(80,10)
\qbezier[1000](0,10)(16,34)(32,10)
\qbezier[1000](16,10)(24,22)(32,10)
\qbezier[1000](0,10)(48,80)(96,10)
\put(-2,0){\small $1$}\put(14,0){\small $2$}\put(30,0){\small $3$}
\put(46,0){\small $4$}\put(62,0){\small $5$}\put(78,0){\small $6$}\put(94,0){\small $7$}
\end{picture}
\end{center}
\caption{A  zigzag stack on $[7]$.}
\label{fig-counter}
\end{figure}

  Lemma \ref{lem-pc} shows that a zigzag stack $S$ can be decomposed into a primary component and a list of zigzag stacks.  Let $S_1, S_2, \ldots, S_k$ be the zigzag stacks on the $C$-intervals $I_1, I_2, \ldots, I_k$. If $S$ is an $m$-reduced zigzag stack, it can be seen that
  all the substructures $S_t$ are $m$-reduced zigzag stacks.
  However, the converse may not be true.  For example, in Figure \ref{fig-counter}, the substructure on  $\{4,5,6\}$ is a $3$-reduced zigzag stack, but the
stack on $\{1,2,3,4,5,6,7\}$ is not a
 $3$-reduced zigzag stack, because for the vertices $3$ and $4$, we have $\ld(3)>0$, $\rd(4)>0$, but $4-3=1<2$,  contradicting the
 condition  $j-i\geq m-1$.

The following lemma gives a necessary and sufficient condition for a zigzag stack $S$ to be $m$-reduced. It shows that the primary
component decomposition of $S$
can be used to restrict the verification of
the following degree conditions to pairs of vertices in each set
$J_t$, where $1\leq t \leq k$:
\begin{enumerate}
  \item[(1)] For $1\le i\le n-m+1$, $\ld(i)+\rd(i+m-1)\leq 2$;
  \item[(2)] For  $1\le i<j\le n$, if $\ld(i)>0$ and  $\rd(j)>0$, then $j-i\geq m-1$.
\end{enumerate}

\begin{lem}\label{lem-mreduced}
  Let $S$ be a zigzag stack on $[n]$, and let $C$ be the primary component of $S$. Let $I_1, I_2,\ldots, I_k$ be the $C$-intervals of $S$, and $J_1, J_2, \ldots, J_k$ be the subsets
  of $[n]$ as given before. Then $S$ is an $m$-reduced zigzag stack if and only if Conditions $(1)$ and $(2)$ hold for each $J_t$ with respect to the degrees in $S$.
\end{lem}

\pf Clearly, we only need to show that if Conditions $(1)$ and $(2)$ hold for each $J_t$, then they hold for any pair of vertices in $[n]$.
The proof consists of two claims. Claim 1. If Condition (1) holds for each $J_t$, then the violation of Condition (1) on $[n]$ implies the violation of Condition (2) on $[n]$.  Claim 2. If Condition (2) holds for each $J_t$, then Condition (2) holds for any pair of vertices in $[n]$. Combining these two claims,
we see that if Conditions (1) and (2) hold for each $J_t$,
then they hold on $[n]$, and hence $S$ is an $m$-reduced zigzag stack.

To prove Claim 1, we suppose that there exist two vertices $i$ and $i+m-1$ such that $\ld(i)+\rd(i+m-1)>2$. Since $\ld(i),\rd(i+m-1)\leq 2$, we see that $\ld(i)>0$ and $\rd(i+m-1)>0$. On the other hand, $i$ and $i+m-1$ cannot be contained in the same interval $J_t$. Otherwise,  Condition (1) would hold for $i$ and $i+m-1$, a contradiction.
Thus there exists at least one vertex $u$ in $C$ such that $i<u<i+m-1$. It follows that either $\ld(u)>0$ or $\rd(u)>0$
 since $C$ is a connected  zigzag stack.
If $\ld(u)>0$, then $(i+m-1)-u<m-1$; If $\rd(u)>0$, then $u-i<m-1$. In either case, Condition (2) is violated.

To verify Claim 2, suppose to the contrary that $i$ and $j$ $(1\le i<j\le n$) are two vertices that violate Condition (2), namely, $\ld(i)>0$, $\rd(j)>0$ but $j-i<m-1$. We further assume that the vertices
 $i$ and $j$ are chosen so that $j-i$ is minimum. Since Condition (2) holds for each $J_t$,  $i$ and $j$ cannot be contained in the same interval $J_t$. Hence there exists at least one vertex $u$ in $C$ such that $i<u<j$.
Noting that $C$ is a connected   zigzag stack, we have either $\ld(u)>0$ or $\rd(u)>0$. If $\ld(u)>0$, then $j-u<j-i$; If $\rd(u)>0$, then $u-i<j-i$. In either case, this contradicts the choice of $i$ and $j$.
So we conclude that any pair of vertices in $S$ satisfies Condition (2). This completes the proof. \qed

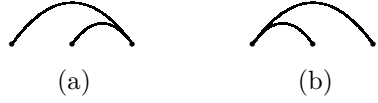
\begin{figure}[H]
\begin{center}
\setlength{\unitlength}{0.5mm}
\begin{picture}(100,20)
\put(0,5){\circle*{1.5}}
\put(16,5){\circle*{1.5}}
\put(32,5){\circle*{1.5}}
\put(64,5){\circle*{1.5}}
\put(80,5){\circle*{1.5}}
\put(96,5){\circle*{1.5}}
\qbezier[1000](16,5)(24,16)(32,5)
\qbezier[1000](0,5)(16,27)(32,5)
\qbezier[1000](64,5)(72,16)(80,5)
\qbezier[1000](64,5)(80,27)(96,5)
\put(12,-7){\small{(a)}}
\put(76,-7){\small{(b)}}
\end{picture}
\end{center}
\caption{Substructures with at least three vertices.}
\label{shape-ab}
\end{figure}

The above lemma shows that an $m$-reduced zigzag stack can be decomposed into a primary component along with a list of $m$-reduced zigzag stacks on the intervals.
We observe that the primary components $C$ has six patterns as shown in Figure \ref{sixcases-1}, where the structure (a) and (b)
 in Figure \ref{shape-ab} stand  for   connected zigzag stacks with
at least three vertices.

\begin{figure}[ht]
\begin{center}
\setlength{\unitlength}{0.5mm}
\begin{picture}(280,80)
\put(3,2){\framebox(10,5)}
\put(19,2){\framebox(10,5)}
\put(35,2){\framebox(10,5)}
\put(51,2){\framebox(10,5)}
\put(67,2){\framebox(10,5)}
\put(97,2){\framebox(10,5)}
\put(113,2){\framebox(10,5)}
\put(129,2){\framebox(10,5)}
\put(145,2){\framebox(10,5)}
\put(161,2){\framebox(10,5)}
\put(191,2){\framebox(10,5)}
\put(207,2){\framebox(10,5)}
\put(223,2){\framebox(10,5)}
\put(239,2){\framebox(10,5)}
\put(255,2){\framebox(10,5)}
\put(271,2){\framebox(10,5)}
\put(11,52){\framebox(10,5)}
\put(27,52){\framebox(10,5)}
\put(43,52){\framebox(10,5)}
\put(59,52){\framebox(10,5)}
\put(104,52){\framebox(10,5)}
\put(120,52){\framebox(10,5)}
\put(136,52){\framebox(10,5)}
\put(152,52){\framebox(10,5)}
\put(206,52){\framebox(10,5)}
\put(222,52){\framebox(10,5)}
\put(238,52){\framebox(10,5)}
\put(254,52){\framebox(10,5)}
\put(0,5){\circle*{1.5}}
\put(16,5){\circle*{1.5}}
\put(32,5){\circle*{1.5}}
\put(48,5){\circle*{1.5}}
\put(64,5){\circle*{1.5}}
\put(94,5){\circle*{1.5}}
\put(110,5){\circle*{1.5}}
\put(126,5){\circle*{1.5}}
\put(142,5){\circle*{1.5}}
\put(158,5){\circle*{1.5}}
\put(188,5){\circle*{1.5}}
\put(204,5){\circle*{1.5}}
\put(220,5){\circle*{1.5}}
\put(236,5){\circle*{1.5}}
\put(252,5){\circle*{1.5}}
\put(268,5){\circle*{1.5}}
\put(8,55){\circle*{1.5}}
\put(24,55){\circle*{1.5}}
\put(40,55){\circle*{1.5}}
\put(56,55){\circle*{1.5}}
\put(101,55){\circle*{1.5}}
\put(117,55){\circle*{1.5}}
\put(133,55){\circle*{1.5}}
\put(149,55){\circle*{1.5}}
\put(203,55){\circle*{1.5}}
\put(219,55){\circle*{1.5}}
\put(235,55){\circle*{1.5}}
\put(251,55){\circle*{1.5}}
\qbezier[1000](0,5)(8,16)(16,5)
\qbezier[1000](32,5)(40,16)(48,5)
\qbezier[1000](32,5)(48,27)(64,5)
\qbezier[1000](0,5)(32,49)(64,5)
\qbezier[1000](94,5)(110,27)(126,5)
\qbezier[1000](110,5)(118,16)(126,5)
\qbezier[1000](142,5)(150,16)(158,5)
\qbezier[1000](94,5)(126,49)(158,5)
\qbezier[1000](204,5)(212,16)(220,5)
\qbezier[1000](188,5)(204,27)(220,5)
\qbezier[1000](236,5)(244,16)(252,5)
\qbezier[1000](236,5)(252,27)(268,5)
\qbezier[1000](188,5)(228,55)(268,5)
\qbezier[1000](24,55)(32,66)(40,55)
\qbezier[1000](24,55)(40,77)(56,55)
\qbezier[1000](8,55)(32,90)(56,55)
\qbezier[1000](117,55)(125,66)(133,55)
\qbezier[1000](101,55)(117,77)(133,55)
\qbezier[1000](101,55)(125,90)(149,55)
\qbezier[1000](203,55)(211,66)(219,55)
\qbezier[1000](235,55)(243,66)(251,55)
\qbezier[1000](203,55)(227,90)(251,55)
\put(33,-5){\tiny ($4$)}
\put(127,-5){\tiny ($5$)}
\put(229,-5){\tiny ($6$)}
\put(33,45){\tiny ($1$)}
\put(126,45){\tiny ($2$)}
\put(228,45){\tiny ($3$)}
\put(6,5){\tiny{$_{T_1}$}}
\put(22,5){\tiny{$_{T_5}$}}
\put(38,5){\tiny{$_{T_1}$}}
\put(54,5){\tiny{$_{T_2}$}}
\put(70,5){\tiny{$_{T_4^*}$}}

\put(100,5.5){\tiny{$_{T_2'}$}}
\put(116,5){\tiny{$_{T_1}$}}
\put(132,5){\tiny{$_{T_5}$}}
\put(148,5){\tiny{$_{T_1}$}}
\put(164,5.5){\tiny{$_{T_4^*}$}}

\put(194,5.5){\tiny{$_{T_2'}$}}
\put(210,5){\tiny{$_{T_1}$}}
\put(226,5){\tiny{$_{T_6}$}}
\put(242,5){\tiny{$_{T_1}$}}
\put(258,5){\tiny{$_{T_2}$}}
\put(274,5.5){\tiny{$_{T_4^*}$}}

\put(14,55.5){\tiny{$_{T_4'}$}}
\put(30,55){\tiny{$_{T_1}$}}
\put(46,55){\tiny{$_{T_2}$}}
\put(62,55.5){\tiny{$_{T_4^*}$}}

\put(107,55.5){\tiny{$_{T_2'}$}}
\put(123,55){\tiny{$_{T_1}$}}
\put(139,55){\tiny{$_{T_4}$}}
\put(155,55.5){\tiny{$_{T_2^*}$}}

\put(209,55){\tiny{$_{T_1}$}}
\put(225,55){\tiny{$_{T_3}$}}
\put(241,55){\tiny{$_{T_1}$}}
\put(257,55.5){\tiny{$_{T_4^*}$}}
\end{picture}
\end{center}
\caption{The primary component decompositions of $m$-reduced zigzag stacks.}
\label{sixcases-1}
\end{figure}

\begin{figure}[h]
  \begin{center}
    \includegraphics[scale=0.5]{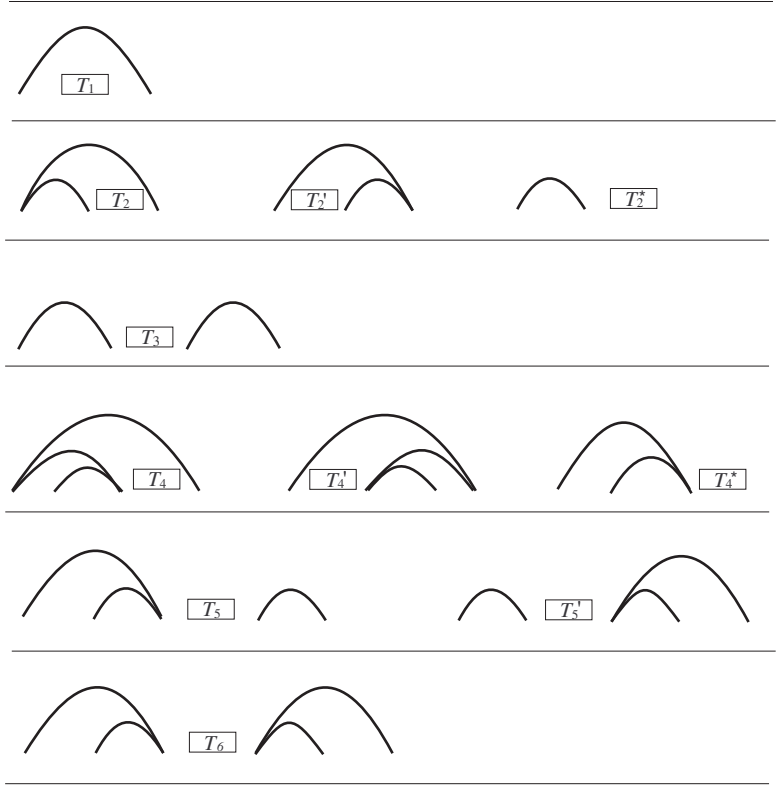}
  \end{center}
  \caption{Types of intervals.}\label{tab-subs}
\end{figure}

We further classify the substructures on the intervals split by $C$ into six classes, see Figure \ref{tab-subs}.
The types of intervals lead us to define types of $m$-reduced
zigzag stacks. If $S$ is a zigzag stack on an interval
of type $T$, then we say that $S$ is also of type $T$.

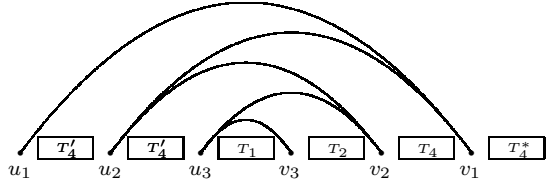
\begin{figure}[ht]
\begin{center}
\setlength{\unitlength}{0.4mm}
\begin{picture}(160,60)
\multiput(0,0)(30,0){6}{\circle*{2}}
\multiput(6,-2)(30,0){6}{\framebox(18,7)}
\qbezier[1000](0,0)(75,100)(150,0)
\qbezier[1000](30,0)(90,80)(150,0)
\qbezier[1000](30,0)(75,60)(120,0)
\qbezier[1000](60,0)(90,40)(120,0)
\qbezier[1000](60,0)(75,22)(90,0)
\put(12,2.5){\tiny{$_{T_4'}$}}
\put(42,2.5){\tiny{$_{T_4'}$}}
\put(72,1){\tiny{$_{T_1}$}}
\put(102,1){\tiny{$_{T_2}$}}
\put(132,1){\tiny{$_{T_4}$}}
\put(162,2){\tiny{$_{T_4^*}$}}
\put(12,2.5){\tiny{$_{T_4'}$}}
\put(42,2.5){\tiny{$_{T_4'}$}}

\put(-4,-6){$_{u_1}$}
\put(26,-6){$_{u_2}$}
\put(56,-6){$_{u_3}$}
\put(86,-6){$_{v_3}$}
\put(116,-6){$_{v_2}$}
\put(146,-6){$_{v_1}$}
\end{picture}
\end{center}
\caption{Decomposition of $m$-reduced zigzag stacks with $\deg(1)=1$.}
\label{typeG}
\end{figure}

Let $S$ be an $m$-reduced zigzag stack on $[n]$, and let $C$ be the primary component of $S$. For given $C$, the types of the intervals created by $C$ are determined.
Assume that $C$ has $k$ vertices, that create $k$ intervals which
are allowed to be empty. In particular, if $C$ containing the vertex $n$, then the $k$-th  interval becomes empty. Let us discuss the case   as shown in Figure \ref{typeG} in detail. The other cases are similar.

When $k=2$, it is easily seen that the two intervals are of types $T_1$ and $T_2^*$ respectively. When $k\ge 3$, the types of the substructures are described below.

\begin{lem}\label{lem-cinterval} Let $S$ be an $m$-reduced
zigzag stack with $\deg(1)=1$, and let $C$ denote the primary component of $S$.
Assume that $C$ contains $k$ vertices with $k\geq 3$.
Suppose that the arcs of $C$ are denoted by $e_1=(u_1, v_1)$, $e_2=(u_2,v_1)$, $e_3=(u_2,v_{2})$, $\ldots$, $e_{k-1}=\big(u_{[\frac{k+1}{2}]}, v_{[\frac{k}{2}]}\big)$.
Then the substructures on the intervals can be described as follows:
\begin{enumerate}
    \item[{\rm (1)}] If $k$ is even, then the substructure on the interval $\langle u_{\frac{k}{2}}, v_{\frac{k}{2}} \rangle$ is of type $T_1$, the one on the interval $\langle v_{\frac{k}{2}}, v_{\frac{k}{2}-1}\rangle$ is of type $T_2$, the ones on the intervals  $\langle u_i,u_{i+1}\rangle$ for $1\le i \le \frac{k}{2}-1$ are of type $T_4'$, the ones on the intervals $\langle v_j,v_{j-1} \rangle$ for $2\le j \le \frac{k}{2}-1$ are of type $T_4$, and the one on the interval after the vertex $v_1$ is of type $T_4^*$.
  \item[{\rm (2)}] If $k$ is odd, then the substructure on the interval $\langle u_{\frac{k+1}{2}}, v_{\frac{k-1}{2}}\rangle$ is of type $T_1$, the one on the interval $\langle u_{\frac{k-1}{2}}, u_{\frac{k+1}{2}}\rangle$ is of type $T_2'$, the ones on the intervals  $\langle u_i,u_{i+1}\rangle$ for $1\le i \le \frac{k-1}{2}-1$ are of type $T_4'$, the ones on the intervals $\langle v_j,v_{j-1}\rangle$ for $2\le j \le \frac{k-1}{2}$ are of type $T_4$, and the one on the interval after the vertex $v_1$ is of type $T_4^*$.
\end{enumerate}
\end{lem}

\pf We only consider the case when $k$ is even. It is clear that $\ld_S(u_i)=0$ for any $1\le i \le k/2$, and $\rd_S(u_1)=1$, $\rd_S(u_i)=2$ for $2\le i \le k/2$. Similarly,  we have $\ld_S(v_j)=2$ for $1\le j \le k/2-1$, $\deg_S(v_{\frac{k}{2}})=1$ and $\rd_S(v_j)=0$ for  $1\le j \le k/2$. Thus for the interval $\langle u_{\frac{k}{2}}, v_{\frac{k}{2}} \rangle$, we have $\ld_S\big(u_{\frac{k}{2}}\big)=0$ and $\rd_S\big(v_{\frac{k}{2}}\big)=0$. By Lemma \ref{lem-mreduced},
 we see that the substructure of $S$ on this interval is of type $T_1$. In other words, the substructure on $\langle u_{\frac{k}{2}}, v_{\frac{k}{2}} \rangle$ can be any $m$-reduced zigzag stack.

For the interval  $\langle v_{\frac{k}{2}}, v_{\frac{k}{2}-1} \rangle$, we have $\ld_S \big(v_{\frac{k}{2}}\big)=1$ and $\rd_S\big(v_{\frac{k}{2}-1}\big)=0$. Thus this interval is of type $T_2$ and the substructure on this interval is also of type $T_2$.

Similarly, the other types of the substructures also can be determined from the degrees of endpoints of the corresponding intervals. This completes the proof. \qed

By Lemma \ref{lem-mreduced}, we see that if $S$ is the substructure in a $C$-interval $\langle u,v\rangle$, then $S$ is an $m$-reduced zigzag stack restricted by $\ld(u)$ and $\rd(v)$; and if $S$ is the substructure in the interval after the last vertex $w$ of the primary component $C$, then $S$ is an $m$-reduced zigzag stack restricted by $\ld(w)$.  For type $T_2$, since $\rd(v)=0$, then $S$ is determined by $\ld(u)$. Noting that $\ld(u)=\ld(w)=1$ for types $T_2$ and $T_2^*$, this implies that the substructures of types $T_2$ and $T_2^*$ are the same class of $m$-reduced zigzag stacks. Since $\ld(u)=\ld(w)=2$ and $\rd(v)=0$ for types $T_4$ and $T_4^*$, we see  that the substructures of types $T_4$ and $T_4^*$ are also the same class of $m$-reduced zigzag stacks. Moreover, if $S$ is of type $T_i$, then the reflection or reversal  of $S$ is of type $T_i'$ $(i=2,4,5)$.

For $1\leq i \leq 6$, let $t_i(n)$ denote the number of $m$-reduced zigzag stacks of type $T_i$ on $[n]$.
Similarly, for $i=2, 4, 5$, let $t_i'(n)$ denote the number of  $m$-reduced zigzag stacks of type $T_i'$ on $[n]$. Moreover, for $i=2, 4$, let $t_i^*(n)$ denote the number of  $m$-reduced zigzag stacks of type $T_i^*$ on $[n]$. It is clear that $t_i(n)=t_i'(n)$
for $i=2,4, 5$. We have shown that $t_i(n)=t_i^*(n)$ for $i=2,4$.
 So we shall identify the types $T_i'$ and $T_i^*$ as $T_i$ and we will  restrict our attention only to substructures of type
$T_i$ $(1\leq i \leq 6)$. For $1 \leq i \leq 6$,
let
 \[
 T_i(x)=\sum_{n=0}^\infty t_i(n)x^n.
 \]

From the primary component decompositions of $m$-reduced zigzag stacks, we find that $Z_m(x)$ can be expressed in terms of $T_i(x)$.

\begin{thm}\label{thm51} We have
\begin{align}\label{gf-zm}
(1-x)Z_m&=1+\frac{x^2T_1T_2}{1-xT_4}+\frac{x^3T_1T_2^2}{1-xT_4}+x^4T_1^2T_3T_4
+\frac{2x^5T_1^2T_2T_4T_5}{1-xT_4}
+\frac{x^6T_1^2T_2^2T_4T_6}{1-xT_4},
\end{align}
where $Z_m$ and $T_i$ stand for the generating functions $Z_m(x)$ and $T_i(x)$, respectively.
\end{thm}

\pf Recall that $z_m(n)$ denotes the number of  $m$-reduced zigzag stacks on $[n]$.
 For $n=0$ , we set $z_m(0)=1$. For $n\ge 1$, let $u_0(n)$ denote the number of these  $m$-reduced zigzag stacks such that $\deg(1)=0$, and define
  \[
  U_0(x)=\sum_{n=0}^\infty u_0(n)x^n.
   \]
   If $1$ is an isolated vertex in  an $m$-reduced zigzag stack
   $S$, we get an $m$-reduced zigzag stack of length $n-1$ by deleting $1$. Thus $u_0(n)=z_m(n-1)$ and
  \[
 U_{0}(x)=xZ_m(x).
 \]

We next divide $m$-reduced zigzag stacks with $\deg(1)\ge 1$ into
six classes $\mathcal{U}_i(n)$ $(1\le i \le 6)$ according to the patterns of  the primary components, see Figure \ref{sixcases-1}.
For each $1\le i\le 6$, denote the  numbers of the $m$-reduced zigzag stacks in $\mathcal{U}_i(n)$ by $u_i(n)$ and define
\[
U_i(x)=\sum_{n=0}^\infty u_i(n)x^n.
\]
  Let $u_i(n,k)$ denote the number of $m$-reduced zigzag stacks in $\mathcal{U}_i(n)$ such that the primary component $C$ contains $k$ vertices. Assume that $C$ contains $k$ vertices $v_1, v_2, \ldots, v_k$.

 Case (1): $\deg(1)=1$. By Lemma \ref{lem-cas2n}, $C$ is of the form as shown in Figure \ref{sixcases-1}. The primary component $C$ creates $k$ intervals. Notice that if $v_k=n$,
 we still consider the
 empty set after the vertex $n$ as an interval.   By Lemma \ref{lem-cinterval}, we see that among these $k$ intervals , there are one interval of type $T_1$, one interval of type $T_2$, and $k-2$ intervals of type $T_4$. Thus
 \[
  u_1(n,k)=\sum_{d_1+\cdots+d_{k}=n-k} t_{1}(d_1)t_2(d_2)t_4(d_3)t_4(d_4)\cdots t_4(d_k),
 \]
where  $d_i$ are nonnegative integers. Therefore,
\begin{align*}
U_1(x)&=\sum_{n=2}^\infty\sum_{k=2}^n u_1(n,k)x^n\\[5pt]
&=\sum_{n=2}^\infty \sum_{k=2}^n \sum_{d_1+\cdots+d_k=n-k}  t_1(d_1)t_2(d_2)t_4(d_3)\cdots t_4(d_k) x^n\\[5pt]
&=\sum_{k=2}^\infty x^k \sum_{n=0}^\infty x^{n} \sum_{d_1+\cdots+d_k=n}  t_1(d_1)t_2(d_2)t_4(d_3)\cdots t_4(d_k) \\[5pt]
&=\sum_{k=2}^\infty x^k\bigg( \sum_{d_1\ge 0} t_1(d_1)x^{d_1} \sum_{d_2\ge 0} t_2(d_2)x^{d_2} \sum_{d_3\ge 0} t_4(d_3)x^{d_3}\ldots \sum_{d_k\ge 0} t_4(d_k) x^{d_k}\bigg)\\[5pt]
&=\sum_{k=2}^\infty x^k T_1(x) T_2(x) T_4^{k-2}(x)\\[5pt]
&=\frac{x^2T_1(x)T_2(x)}{1-xT_4(x)}.
\end{align*}

Case (2): $\deg(1)=2$ and $\deg(v_k)=1$. In this case,
we have $k\ge 3$ and $C$ creates $k$ intervals, in which there are one interval of type $T_1$, two of type $T_2$ and $k-3$ of type $T_4$. It yields that
\[
 u_2(n,k)=
 \sum_{d_1+\cdots+d_{k}=n-k} t_{1}(d_1)t_2(d_2)t_2(d_3)t_4(d_4)t_4(d_5)\cdots t_4(d_{k}),
 \]
where  $d_i$ ranges over nonnegative integers.  Hence
\[
U_2(x)=\frac{x^3T_1T_2^2}{1-xT_4}.
\]

Case (3): $k=4$, $(1, u)$, $(1,w)$ and $(v, w)$ are arcs of $S$ such that $u<v<w\leq n$ and $\deg(u)=\deg(v)=1$. These four vertices
$1, u, v, w$  create four intervals $\langle 1, u\rangle$, $\langle u, v\rangle$, $\langle v, w\rangle$ and the interval after the vertex $w$.  In these intervals, there are two of type $T_1$, one of type $T_3$ and one of type $T_4$. Thus we get
\[
u_3(n)=
 \sum_{d_1+d_2+d_3+d_4=n-4} t_{1}(d_1)t_1(d_2)t_2(d_3)t_4(d_4),
 \]
where  $d_1, d_2, d_3$ and $d_4$ are nonnegative integers.
It follows that
\[
U_3(x)=x^4T_1^2T_3T_4.
\]

Cases (4) and (5) are symmetric, thus $u_4(n)=u_5(n)$ and $U_4(x)=U_5(x)$. Let us consider Case (4), which means that  $(1, v_2)$, $(1,v_k)$ and $(v_3, v_k)$ are arcs of $S$ such that $1<v_2<v_3<v_k\le n$, $\deg(v_2)=1$ and $\deg(v_3)=2$.
In this case $k\ge 5$ and the vertices of $C$ create $k$ intervals, where the interval
after $v_k$ is allowed to be empty. Among these $k$ intervals, there are two   of type $T_1$, one of type $T_2$, one of type $T_5$ and $k-4$ of type $T_4$. So we have
\[
u_4(n,k)=
 \sum_{d_1+\cdots+d_k=n-k} t_{1}(d_1)t_1(d_2)t_2(d_3)t_5(d_4)t_4(d_5)\cdots t_4(d_k),
 \]
 where $d_1, d_2, \ldots, d_k$ are nonnegative integers.
Therefore,
\[
U_4(x)=\frac{x^5T_1^2T_2T_4T_5}{1-xT_4}.
\]

Case (6): $\deg(1)=2$, $(1,v_i)$, $(1,v_k)$, $(v_{i+1},v_k)$ are arcs of $S$ such that $1<v_i<v_{i+1}<v_k$ and $\deg(v_i)=\deg(v_{i+1})=2$.  In this case $k\ge 6$ and $3\leq i\leq k-3$. When $i$ is given, these $k$ vertices create $k$ intervals. It can be seen that among these $k$ intervals, there are two   of type $T_1$, two of type $T_2$, one of type $T_6$ and $k-5$ of type $T_4$.  Thus,
\[
u_6(n,k)=
 \sum_{d_1+\cdots+d_k=n-k}(k-5)t_{1}(d_1)t_1(d_2)t_2(d_3)t_2(d_4)t_6(d_5)t_4(d_6)\cdots t_4(d_k),
 \]
which implies
\[
U_6(x)=\frac{x^6T_1^2T_2^2T_4T_6}{1-xT_4}.
\]

In summary, we find that
\begin{align}\label{gf-gl}
Z_m(x)&=1+U_{0}(x)+U_{1}(x)+U_{2}(x)+U_{3}(x)+U_{4}(x)+U_{5}(x)+U_{6}(x)
\nonumber\\[5pt]
&=1+xZ_m+\frac{x^2T_1T_2}{1-xT_4}+\frac{x^3T_1T_2^2}{1-xT_4}+x^4T_1^2T_3T_4
+\frac{2x^5T_1^2T_2T_4T_5}{1-xT_4}
+\frac{x^6T_1^2T_2^2T_4T_6}{1-xT_4},
\end{align}
which completes the proof. \qed

\section{Generating functions of substructures}
\label{sec-sub}

In this section, we analyze the substructures of types $T_i$ $(1\le i \le 6)$ and derive their generating functions in terms of $Z_m(x)$.
To do so, we shall consider two classes of $m$-reduced
zigzag stacks, which will be called $m$-reduced zigzag stacks of
types $G$ and $H$.

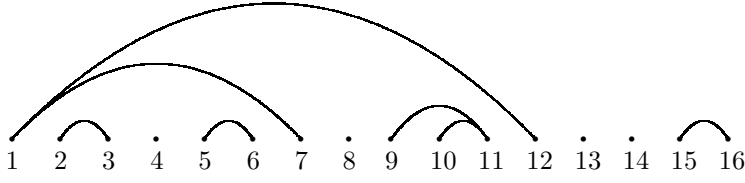
\begin{figure}[h]
\begin{center}
\setlength{\unitlength}{0.4mm}
\begin{picture}(240,70)

\multiput(0,10)(16,0){8}{\circle*{2}}
\qbezier[1000](16,10)(24,22)(32,10)
\qbezier[1000](80,10)(72,22)(64,10)
\qbezier[1000](126,10)(142,32)(158,10)
\qbezier[1000](142,10)(150,22)(158,10)
\qbezier[1000](0,10)(87,100)(174,10)
\qbezier[1000](0,10)(48,60)(96,10)
\qbezier[1000](222,10)(230,22)(238,10)
\put(-2,0){\small $1$}\put(14,0){\small $2$}\put(30,0){\small $3$}
\put(46,0){\small $4$}\put(62,0){\small $5$}\put(78,0){\small $6$}
\put(94,0){\small $7$}\put(110,0){\small $8$}

\multiput(126,10)(16,0){8}{\circle*{2}}
\put(124,0){\small $9$}\put(139,0){\small $10$}\put(155,0){\small $11$}
\put(171,0){\small $12$}\put(187,0){\small $13$}\put(203,0){\small $14$}\put(219,0){\small $15$}\put(235,0){\small $16$}

\end{picture}
\end{center}
\caption{A $3$-reduced zigzag stack.}
\label{fig-GH}
\end{figure}

More precisely, an $m$-reduced zigzag stack on $[n]$ is said to be of type $G$ if $\deg(1)\le 1$, or of type $H$ if $\deg(1)\le 1$ and $\deg(n)\le 1$.
  For example, given the $3$-reduced zigzag stack   in Figure \ref{fig-GH}, the zigzag stack on $\{9, 10, 11\}$ is of type $G$.  The zigzag stack on $\{14,15,16\}$ is also of type $G$. Moreover, the substructure on $\{2,3,4,5,6\}$ is of type $H$.

We find that the substructure on each interval of type $T_i$ $(1\le i \le 6)$  can be characterized as follows.

\begin{thm}\label{thm-subs}
Let $S$ be an $m$-reduced zigzag stack.
Let $I$ be an interval of $S$ of type $T_i$ $(1\leq i \leq 6)$.
Let $T$ be the substructure of $S$ on the interval $I$.  Then $T$ can be described as follows.
\begin{enumerate}
          \item[{\rm (1)}] If $T$ is of type $T_1$, then $T$ can be any $m$-reduced zigzag stack.
          \item[{\rm (2)}] If $T$ is of type $T_2$, then $T$ may be empty, or consists of not more than $m-3$ isolated vertices, or $m-2$ isolated vertices followed by an $m$-reduced zigzag stack of type $G$.
          \item[{\rm (3)}] If $T$ is of type $T_3$, then $T$
            consists of $k$ isolated vertices with $m-2 \leq k \leq 2m-5$, or $T$ is a substructure
             beginning with $m-2$ isolated vertices followed by an $m$-reduced zigzag stack of type $H$ and ending with $m-2$ isolated vertices.
          \item[{\rm (4)}] If $T$ is of type $T_4$, then $T$ is empty, or consists of not more than $m-2$ isolated vertices, or $m-1$ isolated vertices followed by an $m$-reduced zigzag stack.
          \item[{\rm (5)}] If $T$ is of type $T_5$, then $T$ can be $k$ isolated vertices with $m-2\leq k \leq 2m-4$, or $m-1$ isolated vertices followed by an $m$-reduced zigzag stack of type $G$ and $m-2$ isolated vertices.
          \item[{\rm (6)}] If $T$ is of type $T_6$, then $T$ could be $k$ isolated vertices with $m-2\leq k \leq 2m-3$, or $m-1$ isolated vertices followed by an $m$-reduced zigzag stack and $m-1$ isolated vertices.
        \end{enumerate}
\end{thm}

\pf  We assume that  $I¡¡= \langle u,v\rangle$, and denote $I$ by $\{w_1,w_2,\ldots,w_\ell\}$ with $w_1=u+1$, $w_2=u+2$, $\ldots$, $w_\ell=v-1$.

If $T$ is of type $T_1$, then $\ld_S(u)=\rd_S(v)=0$. By Lemma \ref{lem-mreduced}, Conditions (1) and (2)  hold on the interval  $J=\{u,w_1, \ldots,w_\ell, v\}$. For vertices $u$ and $v$, since $\ld_S(u)=\rd_S(v)=0$, there is no restriction with respect to Condition (2).  On the other hand, Condition (1) becomes $\rd_S(u+m-1)\le 2$ and $\ld_S(v-m+1)\le 2$, which are automatically satisfied since $S$ is zigzag.
Hence Conditions (1) and (2) hold on $I=\{w_1, \ldots,w_\ell\}$. In other words,  $T$ can be any $m$-reduced zigzag stack.

If $T$ is of type $T_2$, we have $\ld_S(u)=1$ and $\rd_S(v)=0$. If $\ell< m-2$,   we claim that  $w_1, w_2, \ldots, w_\ell$ are isolated vertices. Otherwise, by Lemma \ref{lem-pc}, there exists an arc $(w_j,w_k)$ in $T$. It follows that $\rd_S(w_j)>0$. Now we have $\ld_S(u)>0$ and $\rd_S(w_j)>0$,  but $w_j-u< m-2$,  which contradicts Condition (2). If $\ell\ge m-2$, by the above claim we see that  $\deg_S(w_i)=0$ for  $1\leq i\leq m-2$. Moreover, by Condition (1), we find that  $\ld_S(u)+\rd_S(w_{m-1})\leq 2$, so that $\rd_S(w_{m-1})\leq 1$. Consequently, the substructure on $\{w_{m-1},w_{m},\ldots,w_{\ell}\}$ is an $m$-reduced zigzag stack of type $G$. Considering all possible structures of $T$ when $T$ is
of type $T_2$, we write
\[
T_2 =\emptyset,\  \bullet,\  \bullet\,\bullet,\  \cdots, \ \underbrace{\bullet\,\bullet\,\cdots\,\bullet}_{m-3}\,,\
\underbrace{\bullet\,\bullet\,\cdots\,\bullet}_{m-2}\,G.
\]

If $T$ is of type $T_3$, we have $\ld_S(u)=\rd_S(v)=1$. If $\ell<m-2$, then $\ld_S(u)>0,\rd_S(v)>0$ but $v-u= \ell+1<m-1$, which violates Condition (2). Thus $\ell$ is at least $m-2$. If $m-2\leq \ell<2m-4$, we claim that $w_1, w_2, \ldots, w_\ell$ are isolated vertices. Otherwise, by Lemma \ref{lem-pc}, there exists an arc $(w_j,w_k)$. It follows that $\rd_S(w_j)>0$ and $\ld_S(w_k)>0$. Since $\ld_S(u)>0$ and $\rd_S(w_j)>0$,  Condition (2) implies that $w_j-u\ge m-1$. Furthermore, since $\ld_S(w_k)>0$ and $\rd_S(v)>0$, we see that $v-w_k\ge m-1$. Hence
\[
\ell=v-u-1>(v-w_k)+(w_j-u)-1\ge 2m-3,
\]
 contradicting the assumption that $\ell < 2m-4$. This proves the claim.

  If $ \ell\ge 2m-4$, by the above claim,  the first $m-2$ vertices $w_1,w_2,\ldots,w_{m-2}$ and  the last  $m-2$ vertices $w_{\ell-m+3}$, $w_{\ell-m+4}$, $\ldots$, $w_{\ell-1}$, $w_\ell$ are isolated vertices. Moreover, according to Condition (1), we have  $\ld_S(u)+\rd_S(w_{m-1})\leq 2$ and  $\ld_S(w_{\ell-m+2})+\rd_S(v)\leq 2$. This yields $\rd_S(w_{m-1})\leq 1$ and $\ld_S(w_{\ell-m+2})\leq 1$,
 so that the substructure  on $\{w_{m-1},w_{m},\ldots,w_{\ell-m+2}\}$ is an $m$-reduced zigzag stack of type $H$. Hence we get
\[
T_3 =\underbrace{\bullet\,\bullet\,\cdots\,\bullet}_{m-2}\,,\ \underbrace{\bullet\,\bullet\,\cdots\,\bullet}_{m-1}\,, \cdots, \ \underbrace{\bullet\,\bullet\,\cdots\,\bullet}_{2m-5}\,,\
\underbrace{\bullet\,\bullet\,\cdots\,\bullet}_{m-2}\,H\,
\underbrace{\bullet\,\bullet\,\cdots\,\bullet}_{m-2}.
\]

Similarly,  if $T$ is of type $T_4$, $T_5$ or $T_6$, all possible
structures of $T$ are given below:
\begin{align*}
T_4& =\emptyset,\  \bullet,\  \bullet\,\bullet,\  \cdots, \ \underbrace{\bullet\,\bullet\,\cdots\,\bullet}_{m-2},\
\underbrace{\bullet\,\bullet\,\cdots\,\bullet}_{m-1}\,Z_m,\\[5pt]
T_5& =\underbrace{\bullet\,\bullet\,\cdots\,\bullet}_{m-1}\,,\ \underbrace{\bullet\,\bullet\,\cdots\,\bullet}_{m}\,, \cdots, \ \underbrace{\bullet\,\bullet\,\cdots\,\bullet}_{2m-4}\,,\
\underbrace{\bullet\,\bullet\,\cdots\,\bullet}_{m-1}\,G\,
\underbrace{\bullet\,\bullet\,\cdots\,\bullet}_{m-2},\\[5pt]
T_6& =\underbrace{\bullet\,\bullet\,\cdots\,\bullet}_{m-1}\,,\ \underbrace{\bullet\,\bullet\,\cdots\,\bullet}_{m}\,, \cdots, \ \underbrace{\bullet\,\bullet\,\cdots\,\bullet}_{2m-3}\,,\
\underbrace{\bullet\,\bullet\,\cdots\,\bullet}_{m-1}\,Z_m\,
\underbrace{\bullet\,\bullet\,\cdots\,\bullet}_{m-1}.
\end{align*}
This completes the proof.  \qed

Let $g(n)$ denote the number of $m$-reduced zigzag stacks of type $G$ on $[n]$, and let $h(n)$ denote the number of
 $m$-reduced zigzag stacks of type $H$ on $[n]$.
 The generating functions of $g(n)$ and $h(n)$ are defined as follows,
\[
G(x)=\sum_{n\ge 0} g(n)x^n, \qquad
H(x)=\sum_{n\ge 0} h(n)x^n.
\]

By Theorem \ref{thm-subs}, we are led to expressions of $T_i(x)$ in terms of $Z_m(x)$, $G(x)$ and $H(x)$.

\begin{thm}\label{lem-gf-ti}
We have
\allowdisplaybreaks
\begin{align}
 T_1(x)&=Z_m(x),\label{eq-int-t-1}\\[8pt]
  T_2(x)&=\frac{1-x^{m-2}}{1-x}+x^{m-2}G(x),\label{eq-int-t-2}\\[8pt] T_3(x)&=\frac{x^{m-2}(1-x^{m-2})}{1-x}+x^{2m-4}H(x), \label{eq-int-t-3}\\[8pt]
   T_4(x)&=\frac{1-x^{m-1}}{1-x}+x^{m-1}Z_m(x),\label{eq-int-t-4}\\[8pt]
  T_5(x)&=\frac{x^{m-1}(1-x^{m-2})}{1-x}+x^{2m-3}G(x),\label{eq-int-t-5}\\[8pt]
T_6(x)&=\frac{x^{m-1}(1-x^{m-1})}{1-x}+x^{2m-2}Z_m(x). \label{eq-int-t-6}
\end{align}
\end{thm}

Next, we consider the primary component decompositions of $m$-reduced zigzag stacks of types $G$ and $H$. We obtain formulas for the generating functions $G(x)$ and $H(x)$ in terms of $Z_m(x)$ and $T_i(x)$  $(1\leq i \leq 6)$.

Recall that $S$ is of type $G$
if $\deg_S(1)\leq 1$.  If  $1$ is an isolated vertex, then $S$ becomes an $m$-reduced zigzag stack of length $n-1$ after deleting the vertex $1$.  If $\deg_S(1)=1$, by Lemma \ref{lem-cas2n}, $C$ is of the form as shown in  Figure \ref{typeG}.
By Lemma \ref{lem-cinterval}, for $k\geq 2$,  $S$ can be decomposed into a connected zigzag stack $C$ and a list of $k$ $m$-reduced zigzag stacks.
Based on this decomposition, we obtain a functional equation
 satisfied by $G(x)$, $Z_m(x)$, $T_2(x)$ and $T_4(x)$.

\begin{lem}\label{lem-G1} We have
\begin{equation}\label{eq-G1G2}
G(x)=1+xZ_m(x)+\frac{x^2Z_m(x)T_2(x)}{1-xT_4(x)}.
\end{equation}
\end{lem}

\pf  Let $S$ be an $m$-reduced zigzag stack of type $G$ on $[n]$, and let $C$ be the primary component of $S$.
Assume that $C$ has $k$ vertices, and
$S$ is decomposed into a connected zigzag stack $C$ and a list of $k$ zigzag stacks $S_1, S_2, \ldots, S_k$, where each $S_i$
is allowed to be empty.  Suppose that $S_i$ has $d_i$ vertices for $1\le i\le k$,  so that we have  $d_1+d_2+ \cdots+d_k=n-k$.

Let $g(n,k)$ denote the number of $m$-reduced zigzag stacks of type $G$ whose  primary component contains $k$ vertices.
Set $g(n,0)=1$.
When $k=1$, by deleting the vertex $1$ from $S$, we obtain an $m$-reduced zigzag stack of length $n-1$. Thus we have $g(n,1)=z_m(n-1)$.

For $k\ge 2$, by Lemma \ref{lem-cinterval}, $C$ creates $k$ intervals, in which there are one of type $T_1$, one of type $T_2$ and $k-2$ of type $T_4$. It follows that
\begin{equation}\label{enus1nk}
g(n,k)= \sum_{d_1+\cdots+d_k=n-k} t_1(d_1)t_2(d_2)t_4(d_3)\cdots t_4(d_k).
\end{equation}
Since
\[
g(n)=\sum_{k=0}^n g(n,k),
\]
 we find that
\allowdisplaybreaks
\begin{align*}
G(x)&=\sum_{n=0}^\infty\sum_{k=0}^n g(n,k)x^n\\[5pt]
&=1+\sum_{n=0}^\infty g(n,1)x^{n}+\sum_{n=0}^\infty \sum_{k=2}^n \sum_{d_1+\cdots+d_k=n-k}  t_1(d_1)t_2(d_2)t_4(d_3)\cdots t_4(d_k) x^n\\[5pt]
&=1+\sum_{n=0}^\infty z_m(n-1)x^{n}+\sum_{k=2}^\infty x^k \sum_{n=0}^\infty x^{n} \sum_{d_1+\cdots+d_k=n}  t_1(d_1)t_2(d_2)t_4(d_3)\cdots t_4(d_k) \\[5pt]
&=1+xZ_m(x)+\frac{x^2T_1(x)T_2(x)}{1-xT_4(x)}.
\end{align*}
Since $T_1(x)=Z_m(x)$, we arrive at the expression (\ref{eq-G1G2})
for $G(x)$.  This completes the proof. \qed

For the $m$-reduced zigzag stacks of type $H$, we obtain the following relation for the generating function $H(x)$.

\begin{lem}\label{lem-G11} We have
\begin{equation}\label{eq-G1G11}
H(x)=1+x+2xG(x)+x^4Z_m^2(x)T_3(x)
+\frac{2x^5Z_m^2(x)T_2(x)T_5(x)}{1-xT_4(x)}
+\frac{x^6Z_m^2(x)T_2^2(x)T_6(x)}{(1-xT_4(x))^2}.
\end{equation}
\end{lem}

\pf  Let $\mathcal{H}(n)$ denote the set of $m$-reduced zigzag stacks of type $H$ on $[n]$. Recall that $h(n)$ denotes the number of $m$-reduced zigzag stack in $\mathcal{H}(n)$. We set $h(0)=1$. It is clear that $h(1)=1$. For $n\geq 2$, according to the degrees of vertices $1$ and $n$, we divide the set of $m$-reduced zigzag stacks of type $H$ into four classes $\mathcal{F}_{1}(n)$, $\mathcal{F}_{2}(n)$, $\mathcal{F}_{3}(n)$ and $\mathcal{F}_{4}(n)$.
For each $1\le i \le 4$, denote the cardinality of $\mathcal{F}_{i}(n)$ by $f_i(n)$ and define the generating function as
\[
F_i(x)=\sum_{n=0}^\infty f_i(n)x^n.
\]

 (1) $\mathcal{F}_{1}(n)=\{S \,|\, S\in  \mathcal{H}(n)\ \mbox{and}\ \deg(1)=0, \deg(n)=0\}$.

Deleting the vertices $1$ and $n$ from a stack  in $\mathcal{F}_{1}(n)$, we obtain an $m$-reduced zigzag stack  on $n-2$ vertices. Conversely, from any given $m$-reduced zigzag stack
 on $n-2$ vertices, we may generate a stack in $\mathcal{F}_1(n)$.
 Thus we have $f_1(n)=z_m(n-2)$ and the generating function
 \[
F_1(x)= x^2Z_m(x).
 \]

 (2) $\mathcal{F}_{2}(n)=\{S \,|\, S\in  \mathcal{H}(n)\ \mbox{and}\ \deg(1)=1, \deg(n)=0\}$.

 By deleting the vertex $n$ from  a stack in $\mathcal{F}_{2}(n)$, we obtain an $m$-reduced zigzag stack of type $G$ of length $n-1$ with $\deg(1)=1$. This process is clearly reversible.  By the definition of type $G$, we see that the number of $m$-reduced zigzag stacks of length $n-1$ with $\deg(1)\le 1$ is $g(n-1)$.
On the other hand, given an $m$-reduced zigzag stack  of length $n-1$ with $\deg(1)= 0$, by deleting the vertex $1$, we are led to an $m$-reduced zigzag stack  of length $n-2$. Thus the number of $m$-reduced zigzag stacks of length $n-1$ with $\deg(1)= 0$ is $z_m(n-2)$. Hence,
 \[
 f_2(n)=g(n-1)-z_m(n-2)
 \]
  and the generating function is
 \[
 F_2(x)=xG(x)-x^2Z_m(x).
 \]

 (3)  $\mathcal{F}_{3}(n)=\{S \,|\, S\in  \mathcal{H}(n)\ \mbox{and}\ \deg(1)=0, \deg(n)=1\}$.

 By reversing the order of vertices $1,2,\ldots,n$, we get a one-to-one correspondence between $\mathcal{F}_{2}(n)$ and $\mathcal{F}_{3}(n)$. Therefore, $f_3(n)=f_2(n)$ and
 \[
F_3(x)=xG(x)-x^2Z_m(x).
 \]

 (4)  $\mathcal{F}_{4}(n)=\{S \,|\, S\in  \mathcal{H}(n)\ \mbox{and}\ \deg(1)=1, \deg(n)=1\}$.

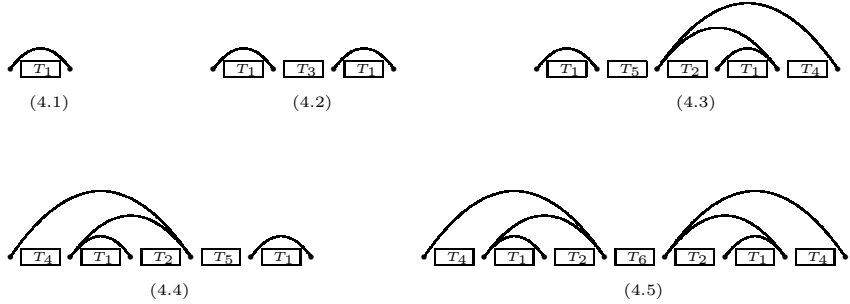
\begin{figure}[ht]
\begin{center}
\setlength{\unitlength}{0.5mm}
\begin{picture}(220,80)
\put(3,3){\framebox(10,4)}
\put(19,3){\framebox(10,4)}
\put(35,3){\framebox(10,4)}
\put(51,3){\framebox(10,4)}
\put(67,3){\framebox(10,4)}
\put(113,3){\framebox(10,4)}
\put(129,3){\framebox(10,4)}
\put(145,3){\framebox(10,4)}
\put(161,3){\framebox(10,4)}
\put(177,3){\framebox(10,4)}
\put(193,3){\framebox(10,4)}
\put(209,3){\framebox(10,4)}
\put(3,53){\framebox(10,4)}
\put(57,53){\framebox(10,4)}
\put(73,53){\framebox(10,4)}
\put(89,53){\framebox(10,4)}
\put(143,53){\framebox(10,4)}
\put(159,53){\framebox(10,4)}
\put(175,53){\framebox(10,4)}
\put(191,53){\framebox(10,4)}
\put(207,53){\framebox(10,4)}
\put(0,5){\circle*{1.5}}
\put(16,5){\circle*{1.5}}
\put(32,5){\circle*{1.5}}
\put(48,5){\circle*{1.5}}
\put(64,5){\circle*{1.5}}
\put(80,5){\circle*{1.5}}
\put(110,5){\circle*{1.5}}
\put(126,5){\circle*{1.5}}
\put(142,5){\circle*{1.5}}
\put(158,5){\circle*{1.5}}
\put(174,5){\circle*{1.5}}
\put(190,5){\circle*{1.5}}
\put(206,5){\circle*{1.5}}
\put(222,5){\circle*{1.5}}
\put(0,55){\circle*{1.5}}
\put(16,55){\circle*{1.5}}
\put(54,55){\circle*{1.5}}
\put(70,55){\circle*{1.5}}
\put(86,55){\circle*{1.5}}
\put(102,55){\circle*{1.5}}
\put(140,55){\circle*{1.5}}
\put(156,55){\circle*{1.5}}
\put(172,55){\circle*{1.5}}
\put(188,55){\circle*{1.5}}
\put(204,55){\circle*{1.5}}
\put(220,55){\circle*{1.5}}
\qbezier[1000](16,5)(24,16)(32,5)
\qbezier[1000](16,5)(32,27)(48,5)
\qbezier[1000](0,5)(24,40)(48,5)
\qbezier[1000](64,5)(72,16)(80,5)
\qbezier[1000](126,5)(134,16)(142,5)
\qbezier[1000](126,5)(142,27)(158,5)
\qbezier[1000](110,5)(134,40)(158,5)
\qbezier[1000](190,5)(198,16)(206,5)
\qbezier[1000](174,5)(190,27)(206,5)
\qbezier[1000](174,5)(198,40)(222,5)
\qbezier[1000](0,55)(8,66)(16,55)
\qbezier[1000](54,55)(62,66)(70,55)
\qbezier[1000](86,55)(94,66)(102,55)
\qbezier[1000](140,55)(148,66)(156,55)
\qbezier[1000](172,55)(196,90)(220,55)
\qbezier[1000](172,55)(188,77)(204,55)
\qbezier[1000](188,55)(196,66)(204,55)
\put(37,-5){\tiny (4.4)}
\put(163,-5){\tiny (4.5)}
\put(5,45){\tiny (4.1)}
\put(75,45){\tiny (4.2)}
\put(177,45){\tiny (4.3)}
\put(6,5){\tiny{$_{T_4}$}}
\put(22,5){\tiny{$_{T_1}$}}
\put(38,5){\tiny{$_{T_2}$}}
\put(54,5){\tiny{$_{T_5}$}}
\put(70,5){\tiny{$_{T_1}$}}
\put(116,5){\tiny{$_{T_4}$}}
\put(132,5){\tiny{$_{T_1}$}}
\put(148,5){\tiny{$_{T_2}$}}
\put(164,5){\tiny{$_{T_6}$}}
\put(180,5){\tiny{$_{T_2}$}}
\put(196,5){\tiny{$_{T_1}$}}
\put(212,5){\tiny{$_{T_4}$}}
\put(6,55){\tiny{$_{T_1}$}}
\put(60,55){\tiny{$_{T_1}$}}
\put(76,55){\tiny{$_{T_3}$}}
\put(92,55){\tiny{$_{T_1}$}}
\put(146,55){\tiny{$_{T_1}$}}
\put(162,55){\tiny{$_{T_5}$}}
\put(178,55){\tiny{$_{T_2}$}}
\put(194,55){\tiny{$_{T_1}$}}
\put(210,55){\tiny{$_{T_4}$}}
\end{picture}
\end{center}
\caption{Five classes of $\mathcal{F}_4(n)$.}\label{apas}
\end{figure}

 In this case, we further divide  $\mathcal{F}_{4}(n)$ into five classes $\mathcal{F}_{4,j}(n)$, where $1\le j\le 5$, see Figure \ref{apas}. For each $1\le j \le 5$, let $f_{4,j}(n)$  be the number of $m$-reduced zigzag stacks in $\mathcal{F}_{4,j}(n)$ and define
 \[
 F_{4,j}(x)=\sum_{n=0}^\infty f_{4,j}(n)x^n.
 \]

 Case (4.1):  $(1,n)$ is an arc of $S$. Then the  interval $ \langle 1,n \rangle$ is of type $T_1$. By Theorem \ref{thm-subs}, the substructure on $ \langle 1,n \rangle$ is an $m$-reduced zigzag stack of length $n-2$. Thus
   $f_{4,1}(n)=z_m(n-2)$ and
   \[
   F_{4,1}(x)=x^2Z_m(x).
   \]

 Case (4.2): $(1, u)$ and $(v, n)$ are two arcs of $S$, where $1< u<v<n$ such that $\deg(u)=\deg(v)=1$. Consider the types of the three intervals. It can be seen that $\langle 1,u\rangle$ is of type $T_1$, $\langle u, v\rangle$ is of type $T_3$, and $\langle v,n\rangle$ is of type $T_1$. Consequently, for $n\geq 4$,
 \[
 f_{4,2}(n)=\sum_{d_1+d_2+d_3=n-4} t_1(d_1)t_1(d_2)t_3(d_3),
 \]
 where $d_1,d_2$ and $d_3$ are nonnegative integers. So we obtain
 that
 \[
 F_{4,2}(x)=x^4 Z_m^2(x) T_3(x).
 \]

 The Cases (4.3) and (4.4) are symmetric. So we have $f_{4,3}(n)=f_{4,4}(n)$ and $F_{4,3}(x)=F_{4,4}(x)$. Let us  consider Case (4.3). In this case, the primary component contains exactly two vertices.
 Suppose that there are $k-2$ vertices
 in the connected component containing $n$, so that there exist $k$
 vertices in these two connected components, where $k\geq 5$. Ignoring the empty interval
 after the vertex $n$, there are $k-1$ intervals created by these
  $k$ vertices. Using the previous arguments, we see that
   the first interval is of type $T_1$, the second is of type $T_5$, and among the other $k-3$ intervals, there
     are one interval of type $T_1$, one   of type $T_2$, and $k-5$ of type $T_4$. Therefore, we get
 \[
 f_{4,3}(n) =\sum_{k=5}^n\sum_{d_1+\cdots+d_{k-1}=n-k} t_1(d_1)t_1(d_2)t_2(d_3)t_5(d_4)t_4(d_5)t_4(d_6)\cdots t_4(d_{k-1}),
 \]
 where $d_1, d_2, \ldots, d_{k-1}$ are nonnegative integers.
  It follows that
 \[
 F_{4,3}(x)=F_{4,4}(x)=\frac{x^5Z_m^2(x)T_2(x)T_5(x)}{1-xT_4(x)}.
 \]

 Case (4.5): Both  the primary component and the
   connected component containing $n$ have at least three vertices. Suppose that these two connected components have
   a total number of $k$ vertices, where $k\ge 6$.
   Given $k$, the primary connected component may have
   $i$ vertices, where $3\leq i \leq k-3$. When $i$ is also given,
   these $k$ vertices  create $k-1$ intervals, not mentioning the
    empty interval after the vertex $n$.
   It can be seen that among these $k-1$ intervals,
   there are two of type $T_1$, two of type $T_2$, $k-6$ of type $T_4$, and one of type $T_6$. Hence
 \[
 f_{4,5}(n)=\sum_{k=6}^n\sum_{d_1+\cdots+d_{k-1}=n-k} (k-5)t_1(d_1)t_1(d_2)t_2(d_3)t_2(d_4)t_6(d_5)t_4(d_6)t_4(d_7)\cdots t_4(d_{k-1}),
 \]
which gives
 \[
 F_{4,5}(x)=\frac{x^6Z_m^2(x)T_2^2(x)T_6(x)}{(1-xT_4(x))^2}.
 \]
So we find that
 \begin{align*}
F_4(x)&= F_{4,1}(x)+F_{4,2}(x)+\cdots+F_{4,5}(x)\\[8pt]
&=x^2Z_m(x)+x^4Z_m^2(x)T_3(x)+\frac{2x^5Z_m^2(x)T_2(x)T_5(x)}{1-xT_4(x)}+\frac{x^6Z_m^2(x)T_2^2(x)T_6(x)}{(1-xT_4(x))^2}.
 \end{align*}
Finally,
we obtain that
\[
H(x)=1+x+F_1(x)+F_2(x)+F_3(x)+F_4(x),
\]
which leads to relation (\ref{eq-G1G11}), and hence the proof is complete. \qed

Substituting the relations \eqref{eq-G1G2} and \eqref{eq-G1G11} into \eqref{eq-int-t-1}--\eqref{eq-int-t-6}, we find that the generating functions $T_2(x)$, $T_3(x)$ and $T_5(x)$ can be expressed in terms of $Z_m(x)$.

\begin{thm} We have
\begin{align}
T_2(x)&=\frac{\big(1-x^{m-1}+x^{m-1}(1-x)Z_m\big)
\big(1-2x+x^m-x^m(1-x)Z_m\big)}
{(1-x)\big(1-2x+x^m-2x^m(1-x)Z_m\big)},\label{re-T2}\\[8pt]
T_3(x)&=x^{m-2}\Big(x^{4m}(x-1)^3Z_m^4+3x^{3m}(x-1)^2(x^m-2x+1)Z_m^3\nonumber\\[5pt]
&\qquad+x^{2m}(x-1)(3x^{2m}-13x^{m+1}+7x^m+9x^2-5x-1)Z_m^2\nonumber\\[5pt]
&\qquad+x^{m}(x^m-2x+1)(x^{2m}-6x^{m+1}+4x^m+2x^2+2x-3)Z_m\nonumber\\[5pt]
&\qquad+(1-x^m)(x^m-2x+1)^2\Big)\nonumber\\[5pt]
&\qquad\quad \quad \Big/\Big((1-x) (1-x^m Z_m) \big( 1- 2 x +x^m-2 x^m(1-x) Z_m\big)^2\Big), \label{re-T3}
\\[8pt]
T_5(x)&=\frac{x^{m-1}\big(1-x^{m-1}+x^{m-1}(1-x)Z_m\big)\big(1-2x+x^m-x^m(1-x)Z_m\big)}
{(1-x)\big(1-2x+x^m-2x^m(1-x)Z_m\big)}. \label{re-T5}
\end{align}
\end{thm}

We have shown that $T_i(x)$ $(1\le i \le 6)$ can be represented in terms of $Z_m(x)$. These relations will be used in the
next section to derive an equation on $Z_m(x)$.

\section{The generating function of $m$-regular linear stacks}

In this section, we use the relations between
$Z_m(x)$ and $T_i(x)$ $(1\leq i \leq 6)$ to derive
an equation satisfied by $Z_m(x)$.
Then we obtain an equation on the generating function
$R_m(x)$  of the number $r_m(n)$ of $m$-regular
linear stacks on $[n]$.
For given $m$, we can deduce a recurrence relation and an asymptotic formula for $r_m(n)$.
For $m=3,4,5,6$, we give asymptotic formulas for $r_m(n)$.

In Section \ref{sec-mred}, we  find the expression \eqref{gf-zm} of $Z_m(x)$ in terms of $T_i(x)$ $(1\leq i\leq 6)$. Moreover, in Section \ref{sec-sub}, we have shown that the generating functions $T_i(x)$ $(1\leq i\leq 6)$ can be expressed in terms of $Z_m(x)$. Combing these relations, we arrive at the following equation satisfied by $Z_m(x)$.

\begin{thm}\label{thm-mreduced}
We have
\begin{equation}\label{eq-zm}
a_5(x)Z_m^5(x)+a_4(x)Z_m^4(x)+a_3(x)Z_m^3(x)+a_2(x)Z_m^2(x)+a_1(x)Z_m(x)+a_0(x)=0,
\end{equation}
where
\begin{align*}
a_0(x)&=(x-1)(x^m - 2x + 1)^2,\\[5pt]
a_1(x)&= - 2x^{3m+1} + x^{3m}+ 12x^{2m+2}-16x^{2m+1}+7x^{2m}-18x^{m + 3}+ 36x^{m + 2}\\
&\qquad - 28x^{m+1}
+7x^m + 4x^4 - 10x^3 + 12x^2 - 6x + 1,\\[5pt]
a_2(x)&= x^m(2x^{3m+1}-15 x^{2m+2}+14 x^{2m+1}-5 x^{2m}+33x^{m+3}-60x^{m+2}+47x^{m+1} \\ & \qquad
-14x^m-16 x^4+39x^3-45x^2+25x-5),\\[5pt]
a_3(x)&= x^{2m}(x-1)(7 x^{2m+1}-28 x^{m+2}+22 x^{m+1}-8x^m+24 x^3-36 x^2+27 x-8),\\[5pt]
a_4(x)&=x^{3m} (x-1)^2 ( 9 x^{m+1}- 16 x^2+ 11 x-4),\\[5pt]
a_5(x)&=4 x^{4m+1}(x-1)^3.
\end{align*}
\end{thm}

From the bijection between $\mathcal{R}_m(n+m-1)$ and $\mathcal{Z}_m(n)$ as given in Theorem \ref{lem-s2n-as2n}, we see that
\[
R_m(x)=1+x+x^2+\cdots+x^{m-2}+x^{m-1}Z_m(x),
\]
or equivalently,
\begin{equation}\label{eq-gr-fr}
Z_m(x)=\frac{(1-x)R_m(x)-(1-x^{m-1})}{(1-x)x^{m-1}}.
\end{equation}
Substituting the above relation into equation \eqref{eq-zm} on $Z_m(x)$,   we arrive at equation \eqref{gf-r-l} given in Theorem \ref{thm22}.

To conclude this paper, we give some special cases of Theorem \ref{thm22}.
When $m=2$, equation \eqref{gf-r-l} on $R_m(x)$ reduces to M\"{u}ller and Nebel's   equation \eqref{gf-exrna} by replacing $R_2(x)$ with $S(z)+1$. For $m=3, 4, 5, 6$, the asymptotic
formulas for $r_m(n)$ have been given in Introduction.

The values of $r_m(n)$ for $m=3, 4, 5, 6$ and
$1\leq n \leq 12$ are given in Table \ref{tab-rmn}.
\begin{table}[h]
\caption{ $r_m(n)$ for $m=3,4,5,6$.}
\begin{center}
\begin{tabular}{|c|c|c|c|c|c|c|c|c|c|c|c|c|}
  \hline
   $n$   & 1& 2& 3 &4&5&6&7&8&9&10&11&12\\
  \hline
   $r_3(n)$   &1&1& 1& 2& 6 &18&54&162& 491&1509 &4692
  &14729 \\
  \hline
   $r_4(n)$   &1&1&1& 1& 2& 6 &18&52&150& 434&1263 &3699
\\
  \hline
   $r_5(n)$   &1&1&1& 1& 1& 2 &6&18&52& 148&422 &1206\\
  \hline
   $r_6(n)$   &1&1&1& 1& 1& 1 &2&6&18& 52&148&420 \\
  \hline
\end{tabular}\label{tab-rmn}
\end{center}
\end{table}

{\noindent \bf Acknowledgements.} This work was supported by the 973 Project, the PCSIRT Project of
the Ministry of Education, and
the National Science Foundation of China.

\end{document}